\documentclass[a4paper]{article}

\usepackage{a4wide}
\usepackage{amsmath,amsthm,amssymb,amsfonts}
\usepackage[hidelinks]{hyperref}
\usepackage{enumerate}
\usepackage{enumitem}
\usepackage{bm}
\usepackage[usenames]{color}
\usepackage[normalem]{ulem}
\usepackage{algorithm}
\usepackage{graphicx,tikz,calc}
\usepackage{graphbox}
\usepackage{hhline}
\usepackage[font=footnotesize,labelfont=bf]{caption}
\usepackage[skip=0pt]{subcaption}
\usetikzlibrary{decorations.pathreplacing}

\numberwithin{equation}{section}
\newtheorem{theorem}{Property}[section]

\theoremstyle{definition}
\newtheorem{definition}[theorem]{Definition}

\newtheorem{example}[theorem]{Example}

\definecolor{darkgreen}{rgb}{0.0,0.5,0.0}
\definecolor{darkblue}{rgb}{0.0,0.0,0.3}
\definecolor{nicosred}{rgb}{0.65,0.1,0.1}
\definecolor{light-gray}{gray}{0.6}
\definecolor{really-light-gray}{gray}{0.8}

\newcommand{\bc}{\ensuremath{\bm{c}}}

\newcommand{\bx}{\ensuremath{\bm{x}}}

\newcommand{\bw}{\ensuremath{\bm{w}}}

\newcommand{\by}{\ensuremath{\bm{y}}}
\newcommand{\bz}{\ensuremath{\bm{z}}}

\newcommand{\bX}{\ensuremath{\bm{X}}}
\newcommand{\bY}{\ensuremath{\bm{Y}}}
\newcommand{\bZ}{\ensuremath{\bm{Z}}}

\newcommand{\tn}[1]{\textnormal{#1}}
\renewcommand{\d}{\tn{d}}

\newcommand{\todo}[1]{\textcolor{black}{#1}}
\newcommand{\change}[1]{\textcolor{black}{#1}}

\newcommand{\rev}[1]{\textcolor{black}{#1}}

\begin{document}
\title{Laguerre tessellations and polycrystalline microstructures: \\ A fast algorithm for generating grains of given volumes}

\author{D.P.~Bourne$^{1}$, P.J.J.~Kok$^{2}$, S.M.~Roper$^{3}$, W.D.T.~Spanjer$^{2}$ \\
	\small $^{1}$Department of Mathematics, Heriot-Watt University, Edinburgh Campus, Edinburgh, EH14 4AS, UK \\
	\small $^{2}$Tata Steel Research and Development, P.O.~Box 10000, 1970 CA IJmuiden, The Netherlands \\
	\small $^{3}$School of Mathematics and Statistics, University of Glasgow, University Place, Glasgow, G12 8QQ, UK
}

\date{June 22, 2020}
%

\maketitle


\begin{abstract}
We present a fast algorithm for generating Laguerre diagrams with cells of given volumes, which can be used for creating RVEs of polycrystalline materials for computational homogenisation, or for fitting Laguerre diagrams to EBSD or XRD measurements of metals. Given a list of desired cell volumes, we solve a convex optimisation problem to find a Laguerre diagram with cells of these volumes, up to any prescribed tolerance. The algorithm is built on tools from computational geometry and optimal transport theory which, as far as we are aware, \change{have not been applied to microstructure modelling before}.
We illustrate the speed and accuracy of the algorithm by generating RVEs with user-defined volume distributions with up to 20,000 grains in 3D. We can achieve volume percentage errors of less than 1\% in the order of minutes on a standard desktop PC. We also give examples of polydisperse microstructures with bands, clusters and size gradients, and of fitting a Laguerre diagram to 3D EBSD measurements of an IF steel.

\end{abstract}



\section{Introduction}

\vspace{-0.2cm}

\subsection{State of the art}
\vspace{-0.1cm}
\change{Voronoi diagrams and their generalizations are often used to represent the microstructure of polycrystalline metals and foams, e.g.,
\cite{AlpersEtAl2015,AltendorfEtAl2014,BarkerBollerheyHamaekers,DepriesterKubler2019,KokKorver99,Liebscher2015,LeonardiScardiLeoni,LLLFP2011,PetrichEtAl2019,PFCRMMO2019,SBDWKKS2016,SWKKS2018,TeferraRowenhorst,WuCaoFan2005}, with individual Voronoi cells representing grains in metals and pores or bubbles in foams. They can be used to generate complex microstructures quickly using a relatively small number of parameters, and they are often used as representative volume elements (RVEs) for computational homogenisation, e.g., \cite{AlsayednoorHarrisonGuo13, GhoshDimiduk, YadegariTurteltaubSuikerKok14}.
}


\change{
	In this paper we focus on the class of weighted Voronoi diagrams known as \emph{Laguerre diagrams} (or \emph{power diagrams} or \emph{radical Voronoi tessellations}), which provide a more accurate description of the geometry of polycrystalline materials than classical Voronoi diagrams \cite{WuCaoFan2005,LLLFP2011}.
However, Laguerre diagrams share the limitation of Voronoi diagrams that
there is not an explicit relation between their generators and their geometric properties, such as the volumes of their cells.
Consequently an active area of research is to develop algorithms for generating Laguerre diagrams with prescribed geometric properties.
}

\change{
	One popular method for approximately controlling the grain size distribution of Laguerre diagrams is using random close packing of spheres and ellipsoids
	\cite{DepriesterKubler2019,LLLFP2011,PFCRMMO2019,WuZhouWangYang2010}.
	This method is inexact, however, since it is impossible to tessellate Euclidean space with spheres or ellipsoids.
}

\change{Several authors have developed methods for fitting Laguerre diagrams to image data (EBSD, XRD) of polycrystals. Here geometric properties such as grain size, centroid location, and aspect ratio are fitted by minimising a measure of the fitting error using
	deterministic and stochastic optimisation methods, e.g.,  \cite{BarkerBollerheyHamaekers,PetrichEtAl2019,SBDWKKS2016,DuanEtAl, Teferra2015}.
	While optimisation methods can give very accurate results, they can also be computationally expensive. Heuristic methods such as \cite{TeferraRowenhorst} trade off fidelity against speed.
}

\change{
	Recently several authors have generated RVEs with curved boundaries and non-convex cells grains using generalisations of Laguerre diagrams such as generalised balanced power diagrams \cite{AlpersEtAl2015,SWKKS2018,TeferraRowenhorst,SBWPBSJ2016}
	and multilevel Voronoi diagrams \cite{KokKorver99,YadegariTurteltaubSuikerKok14}. Complex geometries can also be created using the open-source software DREAM.3D \cite{Dream3D} and \rev {Neper \cite{Neper}}.
}

\change{
	We will discuss some of these methods in more detail and compare them with ours in the Discussion section.
}

\vspace{-0.2cm}

\subsection{Goal of this paper}
\vspace{-0.1cm}
 The goal of this paper is to develop algorithms for creating Laguerre diagrams with user-defined cell size distributions. Our motivation comes from steel modelling. We wish to generate realistic RVEs of single- and multi-phase steels for computational homogenisation simulations. Unlike much of the literature on Laguerre modelling of polycrystals \cite{AlpersEtAl2015,SDKSJ2017,SWKKS2018,TeferraRowenhorst}, our primary aim is not to fit Laguerre diagrams to EBSD or XRD data, but rather to create a tool for generating a rich family of (possibly never-observed) microstructures, which can be combined with multiscale simulations to optimise grain geometries and lead to the development of new alloys. \todo{Having said that, our algorithms are also very well suited for generating Laguerre diagrams with texture intensities that match EBSD data, as we demonstrate in Example \ref{Ex:DataFitting}.} With these applications to steel in mind, we often refer to Laguerre cells as grains, although our results can be applied more generally to other polycrystalline metals and to foams.

\vspace{-0.2cm}

\subsection{Contributions and outline of the paper}
\vspace{-0.1cm}
In Section \ref{Sec:Laguerre} we recall the definition and some important properties of Laguerre diagrams. In particular Property \ref{thm:2} forms the basis of our work. Section \ref{Sec:MainResults} includes our main result, Algorithm 2, for generating `regularised' Laguerre diagrams with grains of prescribed volumes, up to any given tolerance. We also provide Algorithm 1 that can be used for fitting a Laguerre diagram to EBSD measurements of grain volumes and centroids (Example \ref{Ex:DataFitting}). 
We discuss practical issues about how the algorithms can be implemented in Section \ref{Sec:Implementation}. Section \ref{Sec:Examples} includes some large examples (10,000+ grains) and run time tests in 3D, including examples of RVEs of Interstitial Free (IF) steels.

The theory underlying the algorithms presented in Section \ref{Sec:Laguerre} 
uses results from computational geometry and \emph{optimal transport theory} \cite{LevySchwindt2018,Santambrogio}, a field of mathematics that has recently enormously grown in importance and found applications in a wide range of areas including data science, economics, image processing, partial differential equations and statistics. We believe however that this is its first application in the steel industry.

\vspace{-0.1cm}

\section{Laguerre diagrams}
\label{Sec:Laguerre}
\vspace{-0.2cm}

\subsection{Notation and definitions}
\vspace{-0.1cm}
Let $\Omega \subset \mathbb{R}^d$ be the region occupied by a metal. We consider both the 2- and 3-dimensional cases ($d=2$ and $d=3$). For simplicity we assume that $\Omega$ is a convex polygon if $d=2$ or a convex polyhedron if $d=3$. In principle the algorithms below can be used for non-convex regions with curved boundaries, but they become harder to implement. In all our examples below we take $\Omega$ to be a \change{rectangular} box. If $U$ is a subset of $\Omega$, let $|U|$ denote its area if $d=2$ or its volume if $d=3$.

\begin{definition}[\cite{AurenhammerKleinLee13,OkabeBootsSugiharaChiu}]
\label{Def:Lag}
Let $\bx_1,\ldots,\bx_n$ be distinct points in $\Omega$ and $w_1,\ldots,w_n$ be real numbers (not necessarily positive).
The \emph{Laguerre diagram} or \emph{power diagram} generated by the weighted points $(\bx_1,w_1),\ldots,(\bx_n,w_n)$
is the tessellation $\{ L_i \}_{i=1}^n$ of $\Omega$ defined by
\begin{equation}
L_i = \left\{ \bx \in \Omega : |\bx - \bx_i|^2 - w_i \le |\bx - \bx_j|^2 - w_j \; \forall \, j \in \{1,\ldots,n\} \right\}.
\end{equation}
We refer to the sets $L_i$ as \emph{Laguerre cells} or \emph{grains}.
\end{definition}

Laguerre diagrams have the following basic properties \cite{AurenhammerKleinLee13,OkabeBootsSugiharaChiu}:
\begin{itemize}[leftmargin=*]
\item Laguerre cells are convex polygons if $d=2$ or convex polyhedra if $d=3$.
\vspace{-0.1cm}
\item The Laguerre cells tessellate $\Omega$, which means that $\bigcup_{i=1}^n L_i = \Omega$ and cells can only intersect along their boundaries.
\vspace{-0.1cm}
\item If all the weights are equal, $w_1=w_2=\ldots=w_n$, then the Laguerre diagram is simply a Voronoi diagram.
\vspace{-0.1cm}
\item Adding a constant to all the weights does not affect the Laguerre diagram, i.e., the weighted points $\{ (\bx_i,w_i) \}_{i=1}^n$ and $\{ (\bx_i,w_i+c) \}_{i=1}^n$ generate the same diagram for any $c \in \mathbb{R}$.
\vspace{-0.1cm}
\item A generator $\bx_i$ \rev{needs} not belong to its Laguerre cell $L_i$.
\vspace{-0.1cm}
\item There can be empty Laguerre cells, $L_j = \emptyset$ for some $j$.
\end{itemize}

Now we recall two advanced properties of Laguerre diagrams, Properties \ref{thm:existence of weights} and \ref{thm:2}.
These are the key ingredients for generating RVEs with grains of given sizes (given areas if $d=2$ or given volumes if $d=3$).
Property \ref{thm:existence of weights} states that there always exists a Laguerre diagram with grains of given sizes. Property \ref{thm:2} gives a constructive way of finding one.

\begin{theorem}[{\cite[p.~96, Corollary 6.1]{AurenhammerKleinLee13}}, \cite{Aurenhammer98}]
\label{thm:existence of weights}
Let $\bx_1,\ldots,\bx_n$ be distinct points in $\Omega$. Let $m_1,\ldots,m_n$ be positive numbers with $\sum_{i=1}^n m_i = |\Omega|$. Then there exist weights $w_1,\ldots,w_n$ such that the Laguerre diagram $\{ L_i \}_{i=1}^N$ generated by $(\bx_1,w_1),\ldots,(\bx_n,w_n)$ has cells of size $m_1,\ldots,m_n$:
\begin{equation}
|L_i| = m_i \quad \textrm{ for all } i \in \{ 1 , \ldots , n \}.
\end{equation}
\end{theorem}

The weights $w_i$ in Property \ref{thm:existence of weights} can be computed using the following result:
\begin{theorem}[{\cite[pp.~98-100]{AurenhammerKleinLee13}},
	\cite{Aurenhammer98}, {\cite[Theorem 2]{Merigot2011}}]
\label{thm:2}
Let $\bx_1,\ldots,\bx_n$ be distinct points in $\Omega$. Let $m_1,\ldots,m_n$ be positive numbers with $\sum_{i=1}^n m_i = |\Omega|$.
Define the function $g:\mathbb{R}^n \to \mathbb{R}$ by
\begin{equation}
\label{eq:g}
g(w_1,\ldots,w_n) =\sum_{i=1}^n (m_i - |L_i|)w_i + \sum_{i=1}^n \int_{L_i} |\bx-\bx_i|^2 \, \d \bx
\end{equation}
where $\{ L_i \}_{i=1}^N$ is the Laguerre diagram generated by $(\bx_1,w_1),\ldots,(\bx_n,w_n)$. Then
\begin{enumerate}[label=(\roman*)]
\item The function $g$ is concave.
\item The gradient of $g$ has components
\begin{equation}
\frac{\partial g}{\partial w_i} = m_i - |L_i|.
\end{equation}
for all $i \in \{ 1 , \ldots , n \}$.
\item If $\bw=(w_1,\ldots,w_n)$ is a critical point of $g$, i.e., if $\nabla g(\bw)=\mathbf{0}$, then
the Laguerre diagram $\{ L_i \}_{i=1}^N$ generated by $(\bx_1,w_1),\ldots,(\bx_n,w_n)$ has cells of size $m_1,\ldots,m_n$:
\begin{equation}
|L_i| = m_i \quad \textrm{ for all } i \in \{ 1 , \ldots , n \}.
\end{equation}
\end{enumerate}
\end{theorem}
Property \ref{thm:2} forms the basis of Algorithms 1 and 2. It means that if we want to generate a Laguerre diagram with grains of given sizes, then we just need to find critical points of $g$. Since $g$ is concave, this is equivalent to maximising $g$, or to minimising $-g$, which is a smooth, unconstrained, convex optimisation problem. Fast numerical methods are available for solving this \cite{BoydVandenberghe}.

\vspace{-0.2cm}

\subsection{Controlling the spatial distribution of grains}
\vspace{-0.1cm}
Property \ref{thm:2} not only allows \change{one} to control the size distribution of grains, \change{it also gives some control over the spatial distribution.}
Given positive numbers $m_1,\ldots,m_n$ with $\sum_{i=1}^n m_i = |\Omega|$,
there are infinitely many Laguerre diagrams $\{ L_i \}_{i=1}^n$ such that $|L_i|=m_i$ for all $i \in \{1,\ldots,n\}$. This can be seen from Property \ref{thm:2}\emph{(iii)}; \emph{any} choice of
distinct points  $\bx_1,\ldots,\bx_n$ \change{can give} a Laguerre diagram with grains of size $m_1,\ldots,m_n$.
In Section \ref{Subsec:InitSeeds} we will show how to choose $\bx_1,\ldots,\bx_n$ to control the spatial distribution of the grains.

\vspace{-0.2cm}

\subsection{Connection with optimal transport theory}
\vspace{-0.1cm}
Properties \ref{thm:existence of weights} and \ref{thm:2} can also be stated in the language of semi-discrete optimal transport theory\footnote{Technical remark: It can be shown that evaluating the optimal transport (Wasserstein) distance $W_2(\mathcal{L}^d \llcorner \Omega ,\sum_i m_i \delta_{x_i})$ between the Lebesgue measure and a discrete measure \change{generates a partition of}  $\Omega$ into Laguerre cells of size $m_1,\ldots,m_n$.}; see, e.g., \cite[Sec.~6.4.2]{Santambrogio}, \cite{LevySemiDiscrete2015}, \cite{Merigot2011}. This connection \change{provides} a way of finding critical points of $g$ using fast modern methods from optimal transport theory \cite{KitagawaMerigotThibert2019,LevySemiDiscrete2015}.
\change{We discuss this connection further at the end of Section
\ref{Subsec:EnergyDecreasing}.}


\vspace{-0.1cm}

\section{Main results}

\label{Sec:MainResults}

\vspace{-0.2cm}

\subsection{Statement of the algorithms}
\label{Sec:Algorithms}
\vspace{-0.1cm}
For concreteness we state the algorithms in three dimensions, but they can also be used in two dimensions (\change{by substituting} \emph{volume} with \emph{area} and \emph{polyhedron} with \emph{polygon} wherever they appear in Algorithms 1 and 2).
Our main result is Algorithm 2. First however we consider a simplified version, Algorithm 1, which will help us to understand the importance of the regularisation step in Algorithm 2. Algorithm 1 can also be used for data-driven modelling to fit a Laguerre diagram to EBSD or XRD measurements of grain volumes and centroids (see Example \ref{Ex:DataFitting}). Algorithm 1 is not new and goes back at least as far as \cite{Aurenhammer98}. Our role is simply to bring it to the attention of the microstructure modelling community.

\begin{algorithm}[H]
\label{Algo1}
\smallskip
\textbf{Input:} A convex polyhedron $\Omega$ representing a sample of metal, a list of volumes $m_1,\ldots,m_n$ such that $m_i>0$ and \change{$\sum_{i=1}^n m_i=|\Omega|$}, and a \change{relative} error tolerance $\varepsilon$.

\smallskip
\noindent
\textbf{Output:}
The generators $(\bx_1,w_1),\ldots,(\bx_n,w_n)$ of a Laguerre diagram $\{L_i\}_{i=1}^n$ such that grain $L_i$ has volume $m_i$ up to $\varepsilon$ \change{relative} error, i.e.,
\change{$\frac{||L_i|-m_i|}{m_i} < \varepsilon$}, for all $i \in \{1,\ldots,n\}$.

\smallskip
\noindent
\textbf{Method:} \newline
\emph{Initialisation.} Choose or randomly select $n$ distinct points $\bx_1,\ldots,\bx_n$ in $\Omega$. \newline
\emph{Optimisation step.} Use a numerical optimisation method to find $\bw =(w_1,\ldots,w_n)$ that maximises the function $g$ defined in equation \eqref{eq:g}. Initialise the optimisation method using the initial guess $\bw_{\text{init}}=\mathbf{0}$ and
terminate it 
using the stopping criterion 
\change{$|\nabla g(\bw)| <  \varepsilon \min_j m_j$}.
\caption{\label{algo:Algo1} Generate a Laguerre diagram with grains of given volumes}
\end{algorithm}

\begin{example}[Example of Algorithm 1]
	\label{Ex:1}
	Figure \ref{Fig:Alg1} shows an example of Algorithm 1 implemented in MATLAB with $n=50$ grains in the square domain $\Omega=[0,1]\times[0,1]$. 
		\begin{figure}[H]
		\centering
		\includegraphics[trim={1cm 0 1cm 0},clip,width=0.39\linewidth]{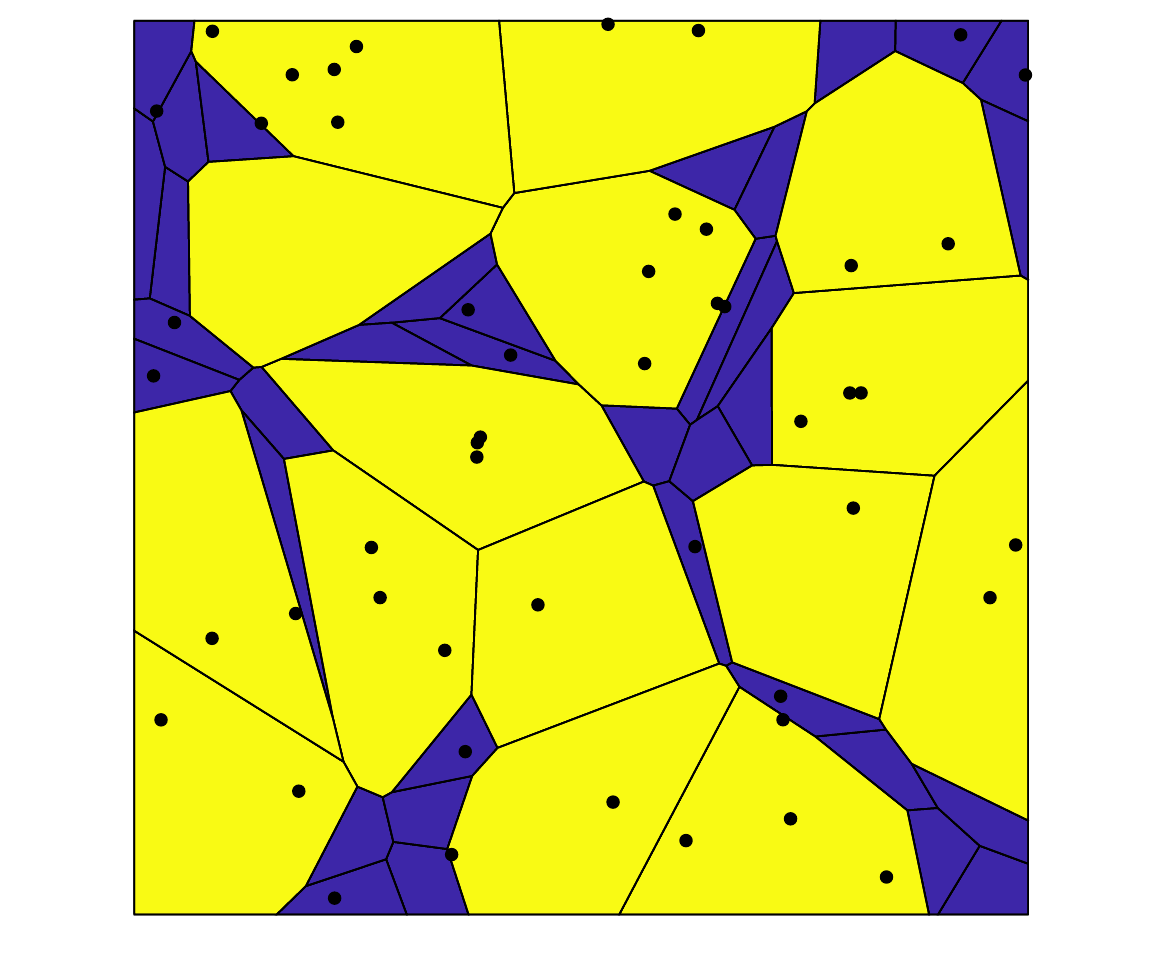}
		\caption{\label{Fig:Alg1} An example of Algorithm 1 with $n=50$ grains in the unit square $\Omega$. There are 35 blue grains and 15 yellow grains. The blue grains have area $x$ and the yellow grains have area $10x$ to within $1\%$ error, where $x=1/185$. The black dots are the locations of the generators $\{ \bx_i \}_{i=1}^{50}$. Notice that not every grain contains its own generator.}
	\end{figure}
	The grains have target areas $m_i=x$ for $i \in \{1,\ldots,35\}$ and $m_i=10x$ for $i \in \{ 36,\ldots,50\}$, where $x=1/185$ so that the total area of all the grains equals the area of $\Omega$. The actual areas of the grains returned by Algorithm 1 are correct to within \change{$1\%$ error ($\varepsilon=0.01$)}. The initialisation step of Algorithm 1 was performed using the MATLAB function \emph{rand} to select $\bx_1,\ldots,\bx_{50}$ (pseudo)randomly from a uniform distribution. While the grains have the correct areas to within $1\%$, the microstructure is somewhat irregular and unrealistic, with some highly elongated grains. This leads us to Algorithm 2, which produces more regular microstructures; compare Figures \ref{Fig:Alg1} and \ref{Fig:Alg2}(i).
\end{example}

\begin{algorithm}[H]
	\label{Algo2}
	\smallskip
	\textbf{Input:} A convex polyhedron $\Omega$ representing a sample of metal, a list of volumes $m_1,\ldots,m_n$ such that $m_i>0$ and \change{$\sum_{i=1}^n m_i=|\Omega|$}, a \change{relative} error tolerance $\varepsilon$, and the number of regularisation steps $K$.
	
	\smallskip
	\noindent
	\textbf{Output:}
	The generators $(\bx_1,w_1),\ldots,(\bx_n,w_n)$ of a regularised Laguerre diagram $\{L_i\}_{i=1}^n$ such that grain $L_i$ has volume $m_i$
	up to $\varepsilon$ \change{relative} error, i.e.,
	\change{$\frac{||L_i|-m_i|}{m_i} < \varepsilon$}, for all $i \in \{1,\ldots,n\}$.
	
	\smallskip
	\noindent
	\textbf{Method:} \newline
	\emph{Initialisation.} Choose or randomly select $n$ distinct points $\bx_1^{(0)},\ldots,\bx_n^{(0)}$ in $\Omega$.
	Initialise the weights to be zero: $\bw^{(0)}=\big(w_1^{(0)},\ldots,w_n^{(0)}\big)=\mathbf{0}$.
	\newline
	\emph{Iteration.} For $k=1,\ldots,K$ do:
	\begin{enumerate}
		\item \emph{Regularisation step.} For $i \in \{1,\ldots,n\}$, define $\bx_i^{(k)}$ to be the centroid of $L_i^{(k-1)}$:
		\begin{equation}
		\label{eq:RegStep}
		\bx_i^{(k)} := \frac{1}{\big|L_i^{(k-1)}\big|} \int_{L_i^{(k-1)}} \bx \, \d \bx
		\end{equation}
		where $\big\{ L_i^{(k-1)} \big\}_{i=1}^n$ is the Laguerre diagram obtained in the previous iteration, which is generated by $\big(\bx_1^{(k-1)},w_1^{(k-1)}\big),\ldots,\big(\bx_n^{(k-1)},w_n^{(k-1)}\big)$.
		\item \emph{Optimisation step.} Use a numerical optimisation method to find $\bw^{(k)} =\big(w_1^{(k)},\ldots,w_n^{(k)}\big)$ that maximises the concave function
		\begin{equation}
		g_k(w_1,\ldots,w_n) =\sum_{i=1}^n \big(m_i - \big|L_i\big|\big)w_i + \sum_{i=1}^n \int_{L_i} \big|\bx-\bx_i^{(k)}\big|^2 \, \d \bx
		\end{equation}
		where $\big\{ L_i \big\}_{i=1}^n$ is the Laguerre diagram generated by $\big(\bx_1^{(k)},w_1\big),\ldots,\big(\bx_n^{(k)},w_n\big)$.
		Initialise the optimisation method using the initial guess $\bw_{\text{init}}=\bw^{(k-1)}$ and
		terminate it using the stopping criterion \change{$\big|\nabla g_k\big(\bw^{(k)}\big)\big| < \varepsilon \min_j m_j$}.
	\end{enumerate}
	\emph{Return.} Output the generators $\big\{ (\bx_i,w_i) \big\}_{i=1}^n = \big\{ \big(\bx_i^{(K)},w_i^{(K)} \big) \big\}_{i=1}^n$.
	\caption{\label{algo:Algo2} Generate a regularised Laguerre diagram with grains of given volumes}
\end{algorithm}

\begin{example}[Example of Algorithm 2]
\label{Ex:2}
Figure \ref{Fig:Alg2} shows an example of Algorithm 2, using $K=100$ iterations, implemented in MATLAB with $n=50$ grains in the square domain $\Omega=[0,1]\times[0,1]$. The grains have target areas $m_i=x$ for $i \in \{1,\ldots,35\}$ and $m_i=10x$ for $i \in \{ 36,\ldots,50\}$, where $x=1/185$ so that the total area of all the grains equals the area of $\Omega$. The actual areas of the grains returned by Algorithm 2 are correct to within \change{$1\%$ error ($\varepsilon=0.01$)}. For the initialisation step we used exactly the same points $\bx_1^{(0)},\ldots,\bx_n^{(0)}$ that \change{were} used for the initialisation step in Example \ref{Ex:1}. Observe from Figure \ref{Fig:Alg2} how the Laguerre diagram becomes more regular as the number of iterations $k$ increases, and how it appears to be converging. The diagram already looks quite regular after just 4 or 5 iterations and the user may be happy to take far fewer than $K=100$ iterations. We discuss how to choose $K$ in the following section.

\begin{figure}[!hbt]
    \centering
    \begin{subfigure}[t]{.32\linewidth}
      \centering
      \subcaption{Iteration $k=1$.}
        \includegraphics[trim={1cm 0 1cm 0},clip,width=\linewidth]{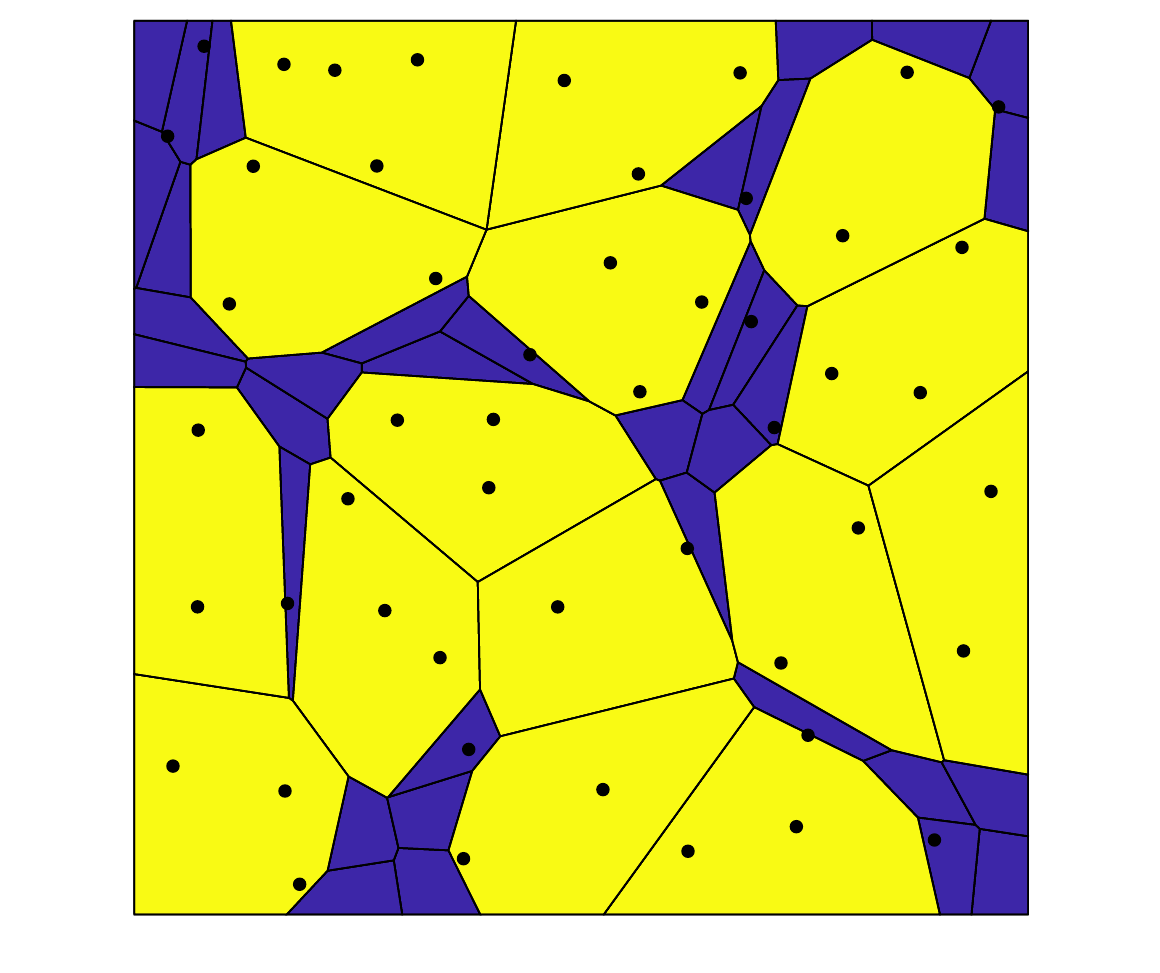}
    \end{subfigure}
    \begin{subfigure}[t]{.32\linewidth}
      \centering
      \subcaption{Iteration $k=2$.}
        \includegraphics[trim={1cm 0 1cm 0},clip,width=\linewidth]{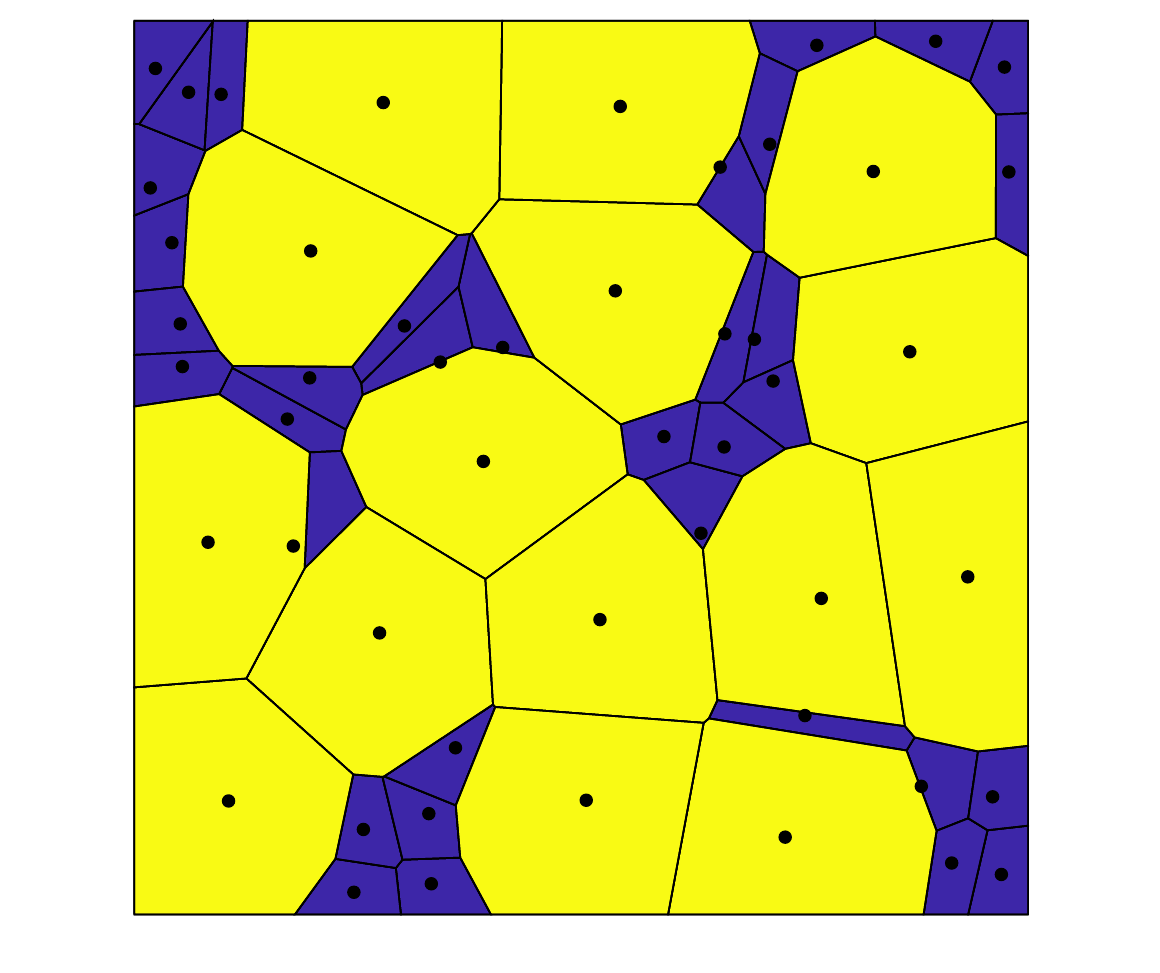}
    \end{subfigure}
    \begin{subfigure}[t]{.32\linewidth}
        \centering
        \subcaption{Iteration $k=3$.}
        \includegraphics[trim={1cm 0 1cm 0},clip,width=\linewidth]{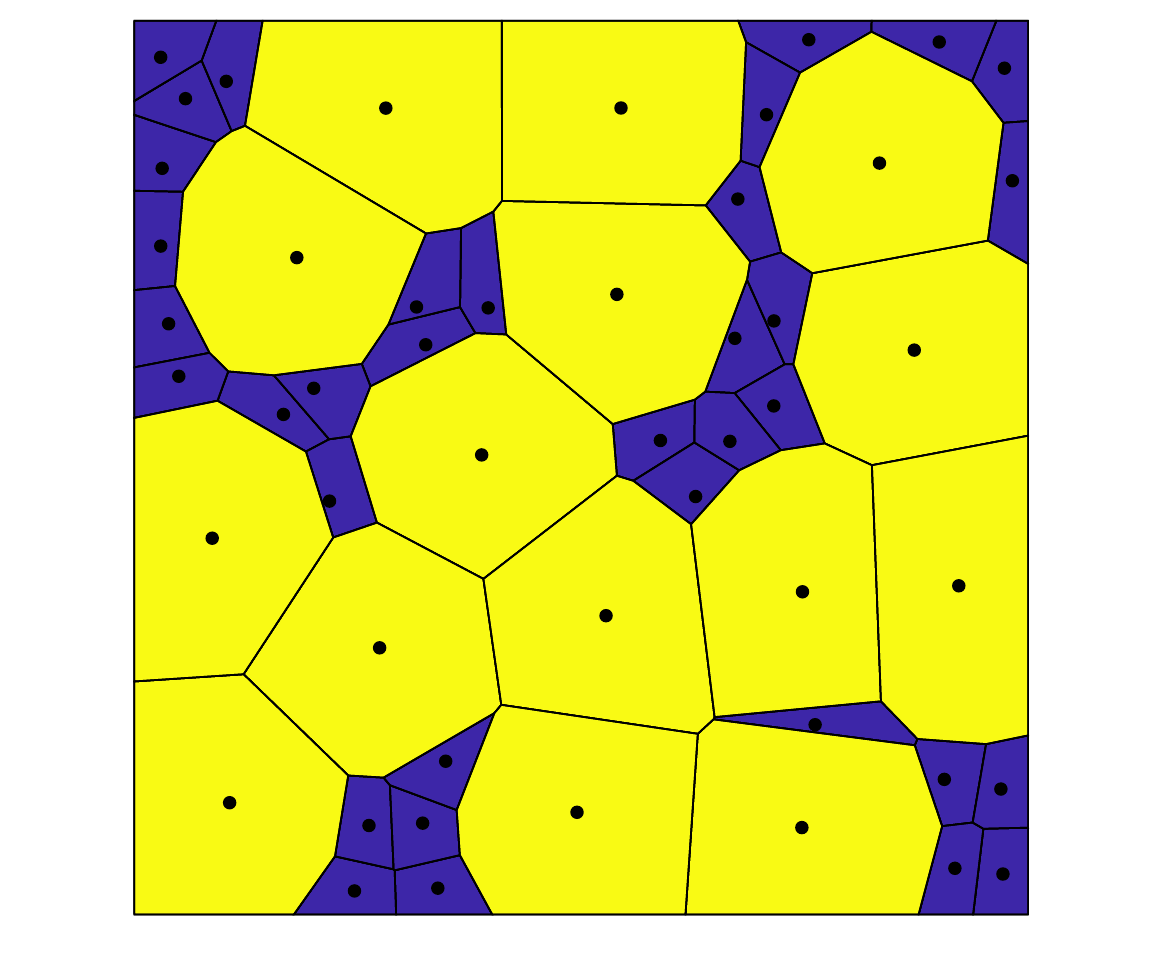}

    \end{subfigure}
    \begin{subfigure}[t]{.32\linewidth}
      \centering
      \subcaption{Iteration $k=4$.}
        \includegraphics[trim={1cm 0 1cm 0},clip,width=\linewidth]{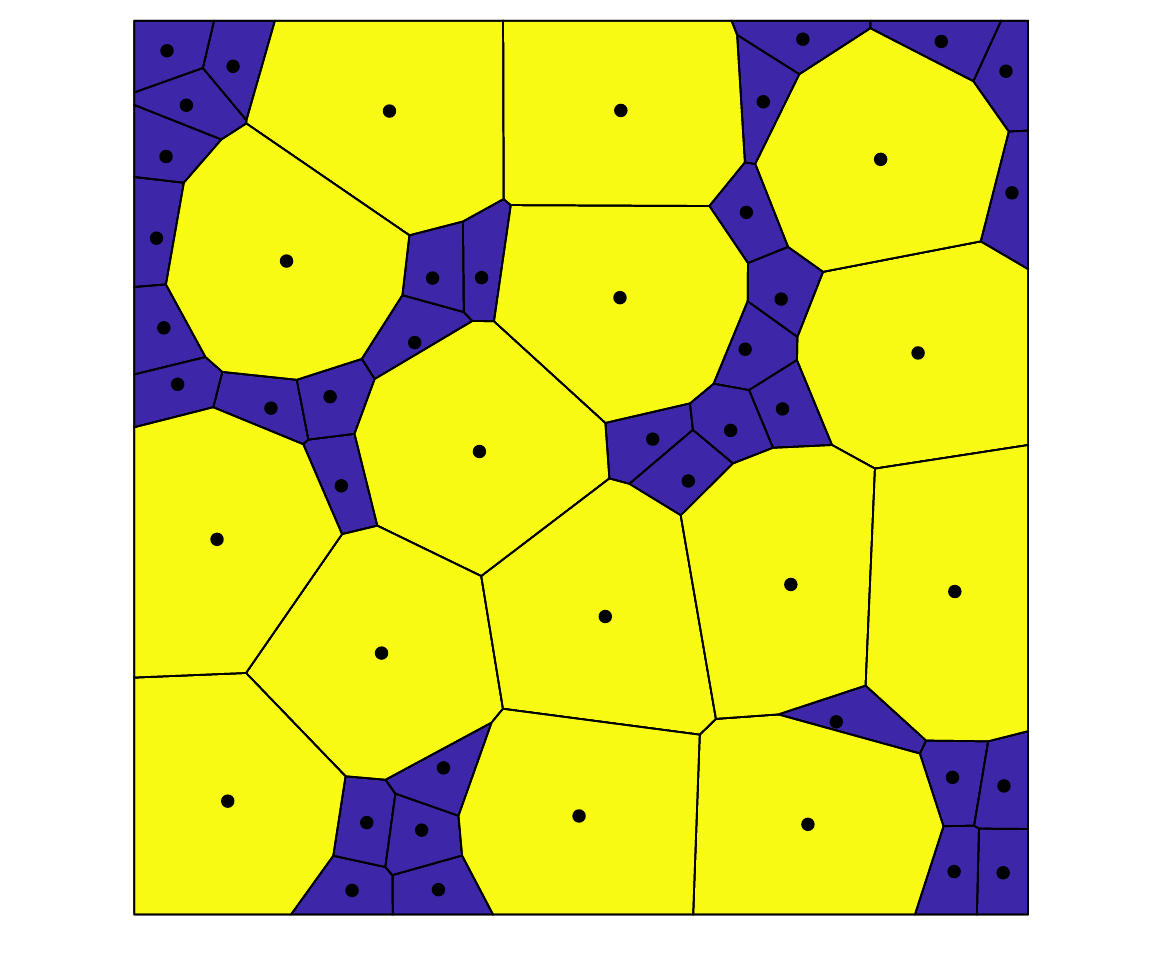}
    \end{subfigure}
    \begin{subfigure}[t]{.32\linewidth}
      \centering
      \subcaption{Iteration $k=5$.}
        \includegraphics[trim={1cm 0 1cm 0},clip,width=\linewidth]{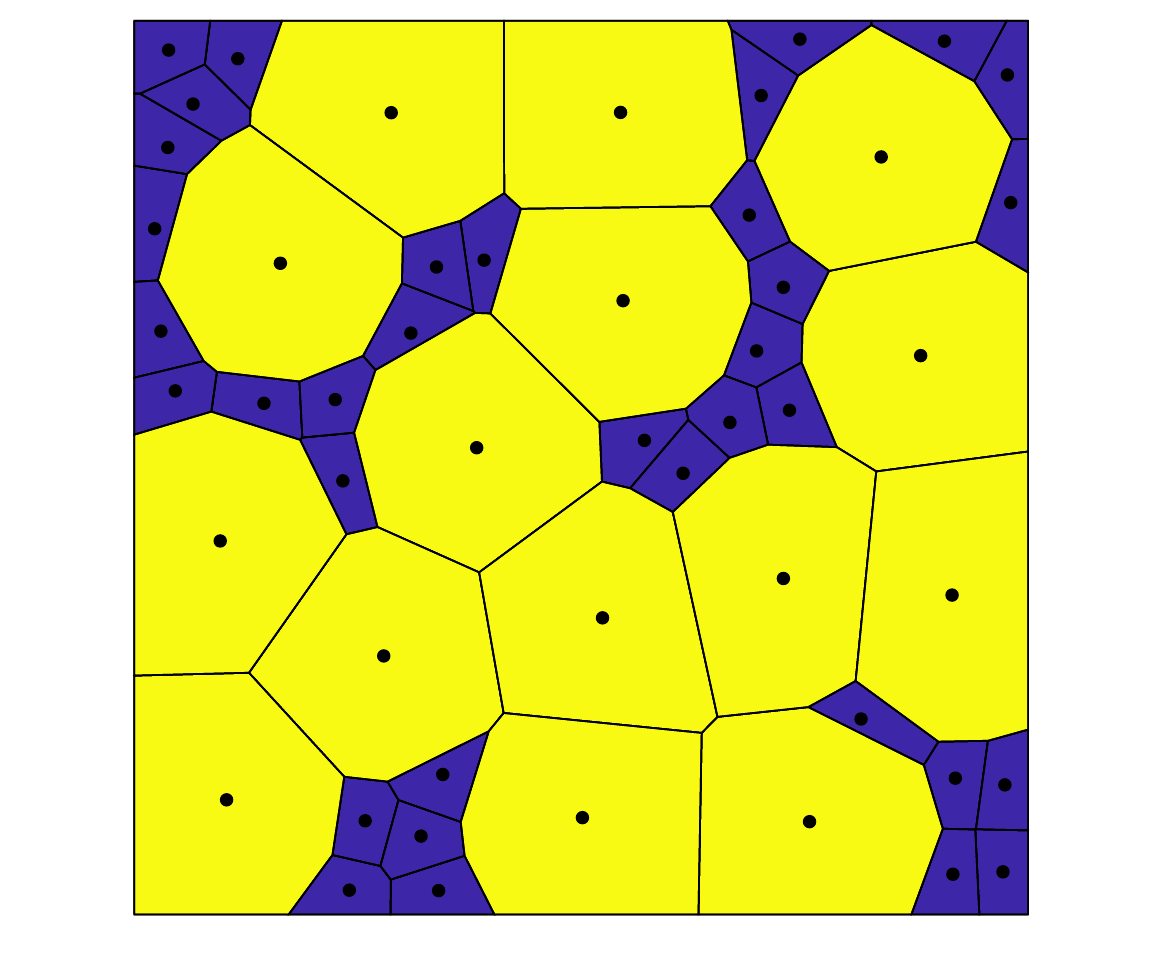}

    \end{subfigure}
    \begin{subfigure}[t]{.32\linewidth}
      \centering
             \subcaption{Iteration $k=10$.}
        \includegraphics[trim={1cm 0 1cm 0},clip,width=\linewidth]{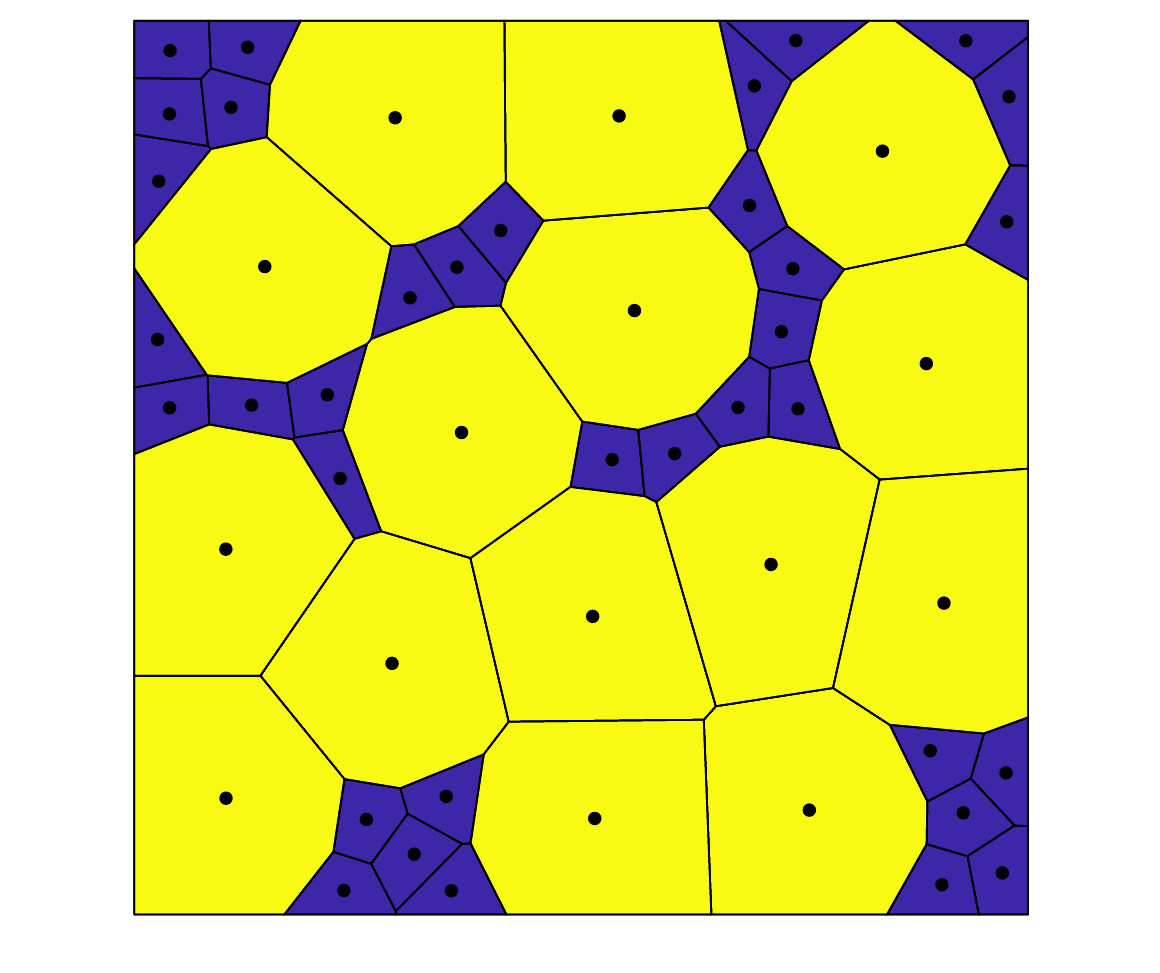}
    \end{subfigure}
    \begin{subfigure}[t]{.32\linewidth}
      \centering
      \subcaption{Iteration $k=25$.}
        \includegraphics[trim={1cm 0 1cm 0},clip,width=\linewidth]{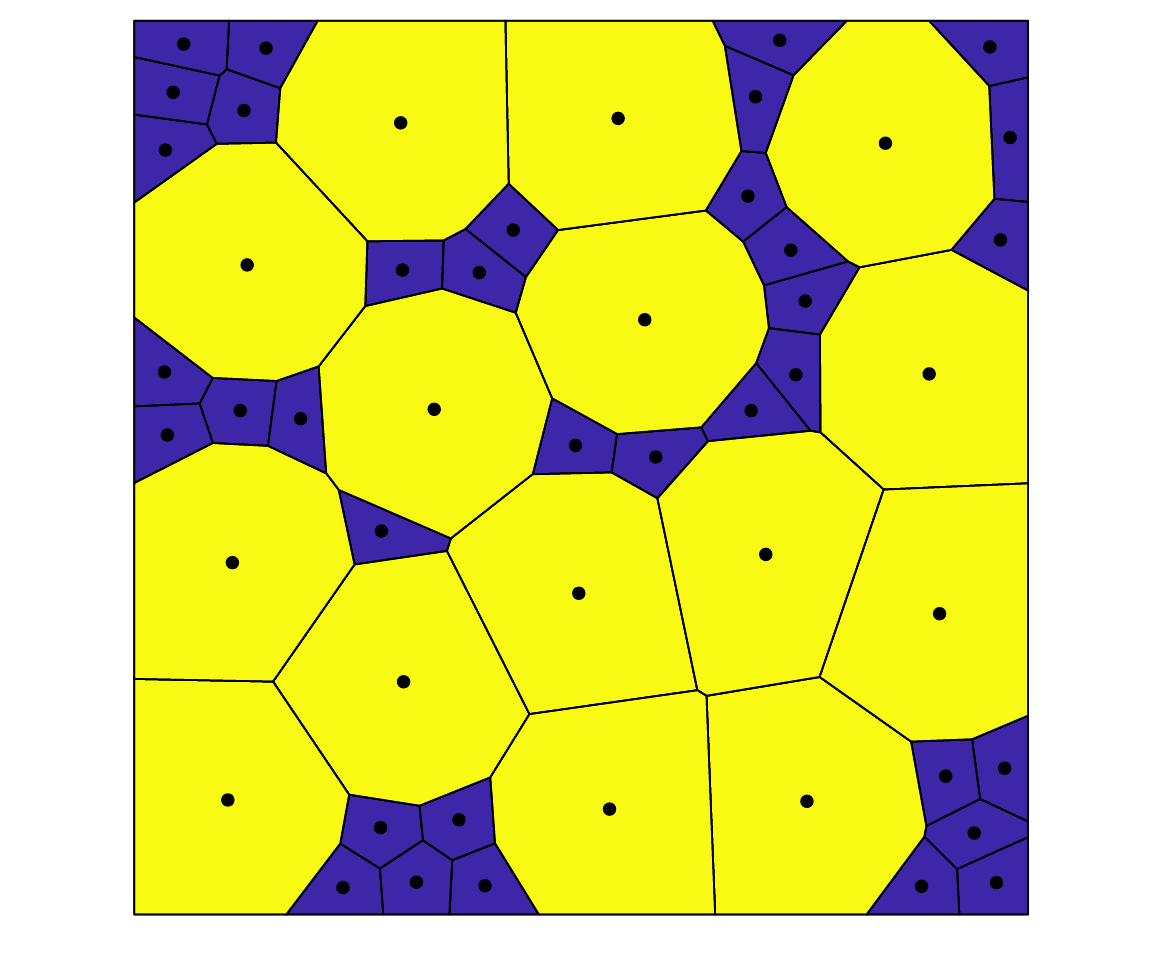}
    \end{subfigure}
    \begin{subfigure}[t]{.32\linewidth}
      \centering
      \subcaption{Iteration $k=50$.}
        \includegraphics[trim={1cm 0 1cm 0},clip,width=\linewidth]{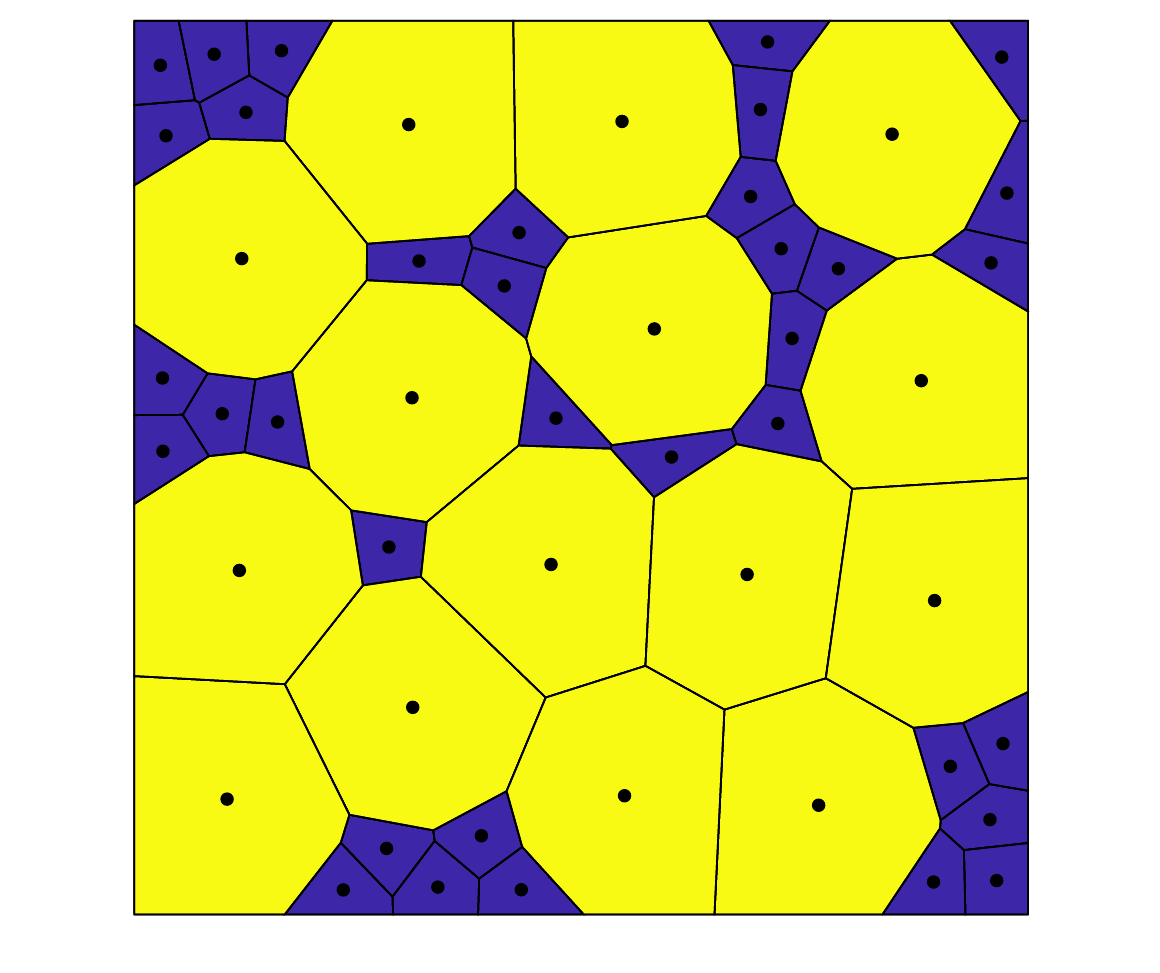}
    \end{subfigure}
    \begin{subfigure}[t]{.32\linewidth}
        \centering
        \subcaption{Iteration $k=100$.}
        \includegraphics[trim={1cm 0 1cm 0},clip,width=\linewidth]{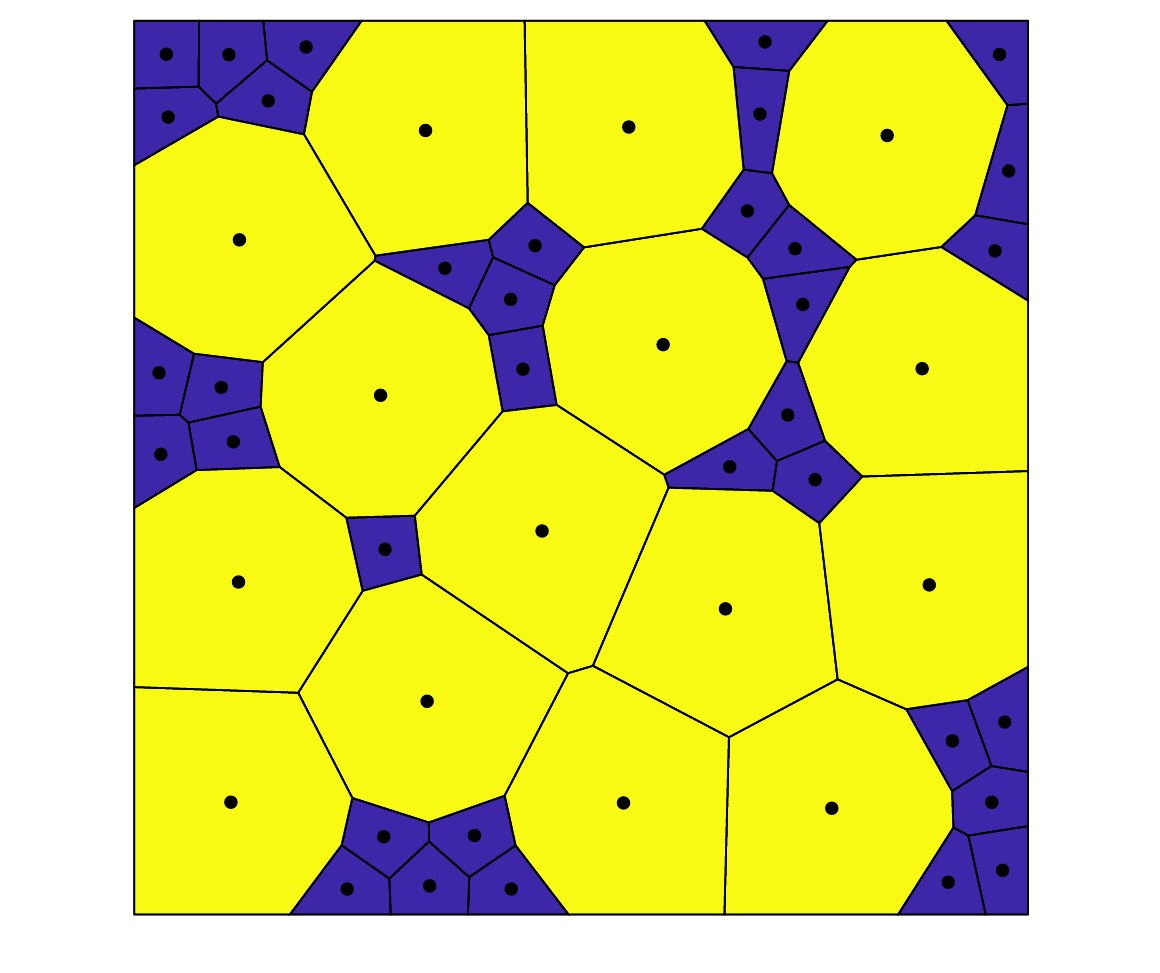}
    \end{subfigure}
    \caption{\label{Fig:Alg2} An example of $K=100$ iterations of Algorithm 2 with $n=50$ grains in the unit square $\Omega$. There are 35 blue grains and 15 yellow grains. The blue grains have area $x$ and the yellow grains have area $10x$ to within $1\%$ error, where $x=1/185$. The black dots are the locations of the generators $\{ \bx_i \}_{i=1}^{50}$. In Figure (i) the generators $\bx_i$ are located at the centroids of their cells $L_i$ to within a distance of $0.002$.}
\end{figure}

\end{example}

\subsubsection{Periodic Laguerre diagrams}

Algorithms 1 and 2 can be modified to create \emph{periodic} Laguerre diagrams for use as RVEs for computational homogenisation (RVEs are usually taken to be periodic to avoid artificial boundary effects). To create periodic diagrams in a rectangular box $\Omega$ of side lengths $l_1,l_2,l_3$, modify Algorithms 1 and 2 as follows. Define the periodic distance between $x,y \in \Omega$ by
\begin{equation}
| x  - y |_\text{per} = \min \{ | x - y + (i l_1, j l_2, k l_3) | : i,j,k \in \mathbb{Z} \}.
\end{equation}
In Algorithms 1 and 2
replace the Laguerre cells $L_i$ by periodic Laguerre cells $\tilde{L}_i$, which are defined by
\begin{equation}
\tilde{L}_i = \left\{ \bx \in \Omega : |\bx - \bx_i|_{\text{per}}^2 - w_i \le |\bx - \bx_j|_{\text{per}}^2 - w_j \; \forall \, j \in \{1,\ldots,n\} \right\}.
\end{equation}
In Algorithm 1 replace $g$ by
\begin{equation}
g_{\text{per}}(w_1,\ldots,w_n) =\sum_{i=1}^n (m_i - |\tilde{L}_i|)w_i + \sum_{i=1}^n \int_{\tilde{L}_i} |\bx-\bx_i|_{\text{per}}^2 \, \d \bx.
\end{equation}
Replace $g_k$ in Algorithm 2 is an analogous way.

\subsection{Properties of the algorithms}

\subsubsection{Convergence of Algorithm 2: centroidal Laguerre diagrams}

\change{In Appendix \ref{Appendix1} we prove that, under a generic assumption,
Algorithm 2 converges as $K \to \infty$. This means that the generator locations $\bx_i^{(k)}$ settle down with increasing iterations}, like we see in Figure \ref{Fig:Alg2}. To be more precise, there exist $\bx_1,\ldots,\bx_n$ such that $\lim_{k \to \infty} \bx_i^{(k)}=\bx_i$ for all $i \in \{1,\ldots,n\}$.
By taking the limit $k \to \infty$ in equation \eqref{eq:RegStep}, we see that
\begin{equation}
\bx_i = \frac{1}{|L_i|} \int_{L_i} \bx \, \d \bx.
\end{equation}
Therefore the generator $\bx_i$ is the centroid of its own Laguerre cell $L_i$ for all $i$.
Such a Laguerre diagram is known as a \emph{centroidal Laguerre diagram} or a \emph{centroidal power diagram}, a term introduced in \cite{BriedenGritzmann};
see also \cite{BourneRoper15,LevyCentroidalPower}.
Centroidal Laguerre diagrams tend to be more regular than non-centroidal Laguerre diagrams, as illustrated by Figure \ref{Fig:Alg2}(i) (centroidal) and
Figure \ref{Fig:Alg1} (non-centroidal).

\subsubsection{Connection with Lloyd's algorithm}

If we omit the optimisation step in Algorithm 2 and set the weights to be zero for all iterations, $\bw^{(k)}=\mathbf{0}$ for all $k$, then we obtain the well-known \emph{Lloyd's algorithm} for computing \emph{centroidal Voronoi tessellations} (Voronoi diagrams where each generator is the centroid of its own Voronoi cell)
\cite{DuFaberGunzburger1999}. Therefore Algorithm 2 can be interpreted as a \emph{generalised Lloyd algorithm with capacity constraints} where cell $L_i$ is constrained to have volume $m_i$. An alternative method for generating centroidal Laguerre diagrams with capacity constraints is given in \cite[Sec.~4]{LevyCentroidalPower}.

\subsubsection{Energy-decreasing property of Algorithm 2}

\label{Subsec:EnergyDecreasing}
Algorithm 2 can also be interpreted as an energy-decreasing optimisation method.
Given $m_1,\ldots,m_n$ with $\sum_i m_i = |\Omega|$, define
\begin{equation}
\label{eq:E0}
E(\bx_1,\ldots,\bx_n) = \min \left\{ \sum_{i=1}^n \int_{U_i} |\bx-\bx_i|^2 \, \d \bx : \{ U_i \}_{i=1}^n \textrm{ is a partition of } |\Omega|, \, |U_i|=m_i \right\}.
\end{equation}
Here the minimum is taken over all possible partitions of $\Omega$, not just Laguerre diagrams. This is an example of an \emph{optimal transport problem}. For example, in two dimensions $E$ could represent the minimum (squared) cost of transporting the recyclable waste generated by a population uniformly distributed over a country $\Omega$ to recycling centres located at $\{ \bx_i \}_{i=1}^n$ with capacities $\{ m_i \}_{i=1}^n$.
It can be shown \cite[Sec.~6.4.1]{AurenhammerKleinLee13} that
\begin{equation}
\label{eq:E}
E(\bx_1,\ldots,\bx_n) = 
\sum_{i=1}^n \int_{L_i} | \bx - \bx_i |^2 \, \d \bx
\end{equation}
where $\{ L_i \}_{i=1}^n$ is the Laguerre diagram with generators $\{ (\bx_i,w_i) \}_{i=1}^n$, where $(w_1,\ldots,w_n)$ is a maximum point of $g$ (defined in \eqref{eq:g}). In other words, $\{ L_i \}_{i=1}^n$ is the solution of the optimal transport problem and all the recyclable waste generated in region $L_i$ should be sent to the recycling centre $\bx_i$.

We could further ask what are the best locations of the recycling centres by considering the optimisation problem
\begin{equation}
\min_{\bx_1,\ldots, \bx_n} E(\bx_1,\ldots,\bx_n).
\end{equation}
This is known as the \emph{optimal location problem} in the economics literature \cite{BollobasStern} and the \emph{quantization problem} in the discrete geometry \cite{Gruber04}, electrical engineering \cite{GrayNeuhoff} and probability literature \cite{GrafLuschgy2000}.
\change{It can be shown that
\begin{equation}
\frac{\partial E}{\partial \bx_i}(\bx_1,\ldots,\bx_n)= 2
\int_{L_i} (\bx_i - \bx) \, \d \bx.
\end{equation}
See for example \cite{LevyCentroidalPower}. Therefore $\nabla E(\bx_1,\ldots,\bx_n)=0$ if and only if $\{\bx_i\}_{i=1}^n$ generate a centroidal Laguerre diagram.}

\change{Thanks to its regularisation step, Algorithm 2 is energy-decreasing in the sense that
 \begin{equation}
 \label{eq:EnergyDecreasing}
 E\big( \bx_1^{(k+1)},\ldots,\bx_n^{(k+1)} \big) \le  E\big( \bx_1^{(k)},\ldots,\bx_n^{(k)} \big).
 \end{equation}
Moreover, under the generic assumption \eqref{eq:ass}, the sequence $( \bx_1^{(k)},\ldots,\bx_n^{(k)} \big)$ converges to a critical point of $E$ (to a local minimum point or saddle point). In other words, it converges to a centroidal Laguerre diagram. We prove these statements in Appendix \ref{Appendix1}. In general\rev{,} Algorithm 2 does not \rev{converge} to a \emph{global} minimum point of $E$ since $E$ is highly non-convex with many critical points; Figure \ref{fig:init_seeds_bimodal} illustrates four different (approximate) critical points of $E$, corresponding to different choices of $( \bx_1^{(0)},\ldots,\bx_n^{(0)})$.}

An alternative method for finding local minima of $E$ is given in \cite[Sec.~4]{LevyCentroidalPower} where, instead of updating $x_i^{(k)}$ using our regularisation step, they update it using a quasi-Newton (L-BFGS) optimisation step applied to $E$.


\section{Implementation}
\label{Sec:Implementation}
In this section we discuss different options for implementing Algorithms 1 and 2, which we \change{did} using MATLAB and Voro++ \cite{Voro++}.

\subsection{Computing Laguerre diagrams}
\label{Subsec:Lag}
One of the main expenses of Algorithms 1 and 2 is the computation of Laguerre diagrams. This happens whenever the objective function $g$ or $g_k$ is evaluated, which could happen many times within a single optimisation step. A Laguerre diagram of $n$ generators can be computed in $\mathcal{O}(n \log n)$ flops in 2D and $\mathcal{O}(n^2)$ flops in 3D \change{\cite[p.~85]{AurenhammerKleinLee13}. (Note that these are \emph{worst-case} optimal run times and in practice the complexity may be better, as we observed in Example \ref{Example:RunTimeTests}. For example, the complexity can be better if each cell has only $\mathcal{O}(n)$ faces instead of the worst-case $\mathcal{O}(n^2)$ \cite[Theorem 6.1]{AurenhammerKleinLee13}.)}
In applications $n$ could be $10,000$ or more, and hundreds or thousands of Laguerre diagrams could be computed in a single run of either algorithm. Therefore it is important to use efficient software.

For our 2D computations we used the MATLAB function \emph{power\_bounded} from the MATLAB File Exchange
\cite{PowerBounded}, which implements Aurenhammer's lifting method \cite{Aurenhammer87} and crops the diagram to a rectangular box $\Omega$.

The function \emph{power\_bounded} is limited to 2D, and so for our 3D computations we used (a slightly modified version of) the C++ library Voro++ \cite{Voro++} combined with a MEX file so that we could run Voro++ via MATLAB. In 3D we also tried the MATLAB function \emph{powerDiagramWrapper} from the MATLAB File Exchange
\cite{PowerDiagramWrapper}, combined with our own code to crop the diagram to a cuboid $\Omega$, but we found Voro++ to be faster. Another advantage of Voro++ is that it can create periodic Laguerre diagrams.

We also used Tata Steel's own in-house Laguerre diagram code to visualise Laguerre diagrams in 2D and 3D.

\subsection{Optimisation methods}
The other main expense of Algorithms 1 and 2 is the optimisation step. For each algorithm this is a smooth, unconstrained, concave maximisation problem and so is very tractable.
We used the MATLAB function \emph{fminunc} to minimise $-g$ and $-g_k$ (and hence maximise $g$ and $g_k$), which uses the \change{BFGS quasi-Newton method by default} \cite{BoydVandenberghe}. This requires an initial guess $\bw_{\text{init}}$ for the minimum point.

\subsubsection{Choice of the initial guess}
\label{Subsec:InitW}
For Algorithm 1 we recommend the initial guess $\bw_{\text{init}}=\mathbf{0}$. 
For data fitting (like Example \ref{Ex:DataFitting}, where the seeds $\bx_i$ are taken from EBSD measurements), if the target grains are relatively spherical, then a better choice may be $(\bw_{\text{init}})_i=m_i/\pi$ in 2D or $(\bw_{\text{init}})_i=(3m_i/(4\pi))^{2/3}$ in 3D. In other words, $(\bw_{\text{init}})_i=r_i^2$ where $r_i$ is the radius of a ball of area $m_i$ in 2D or volume $m_i$ in 3D. This is motivated by sphere-packing methods \cite{DepriesterKubler2019,LLLFP2011,PFCRMMO2019,WuZhouWangYang2010}.

\begin{figure}[tbp]
    \centering
    \includegraphics[width=0.75\textwidth]{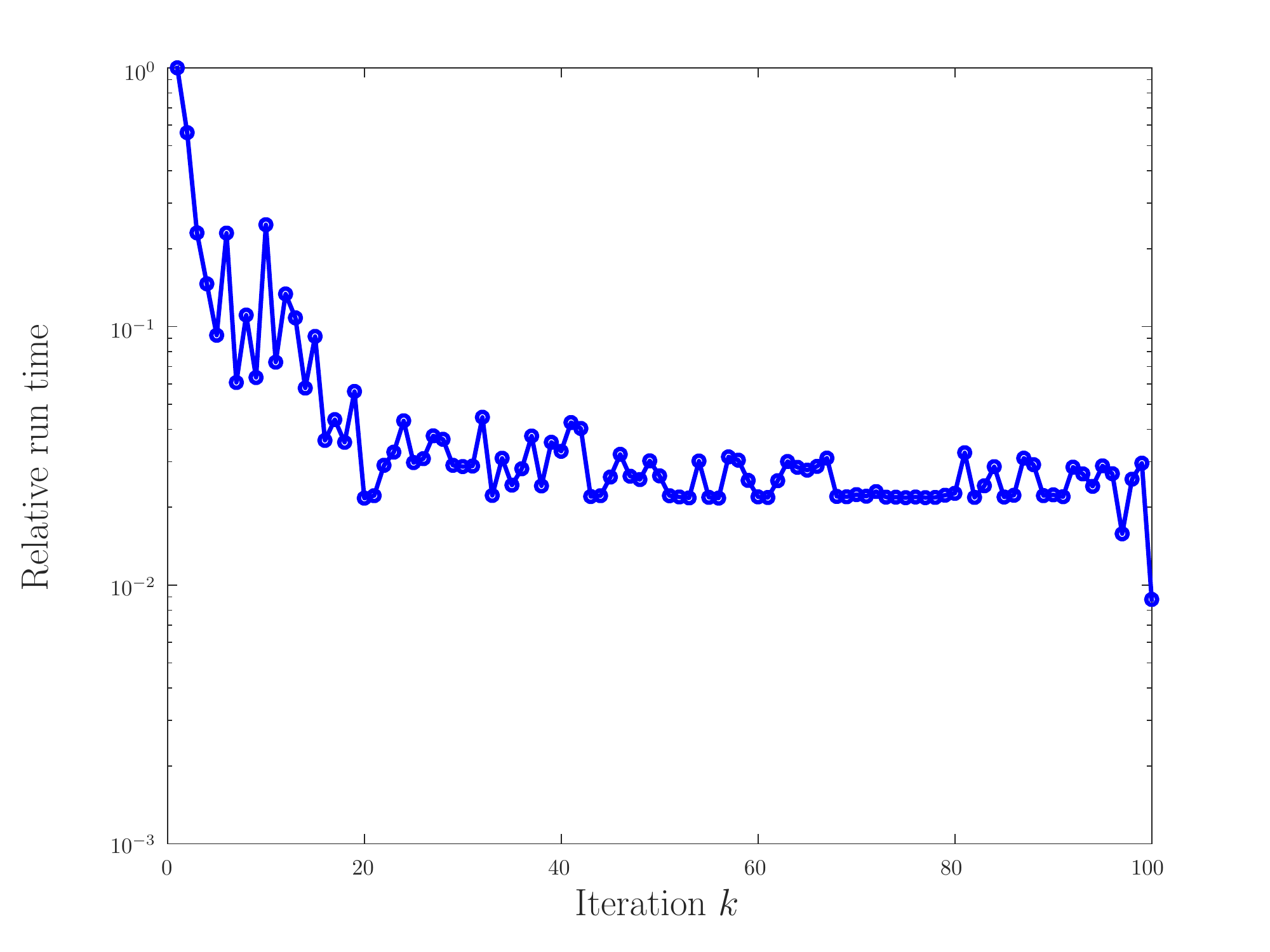}
    \caption{\label{Fig:RelRunTimes} \change{The relative run times for each optimisation step of Algorithm 2. In this example, the parameters are the same as those in Example \ref{Ex:2}, except that the number of grains is $n=500$. The relative proportion of small and large grains is the same as in Example \ref{Ex:2} and the ratio of the size of the largest grain to the size of the smallest grain is $10$. The initial seed locations are randomly distributed in the unit square. The $y$-axis displays $t_k/t_1$, where $t_k$ is the run time for the optimisation step of iteration $k$. After 12 iterations the time per iteration falls to below one-tenth the time of the first iteration, and after 18 iterations it falls to below one-twentieth. This is because the initial guess for the weights in the inner iteration is taken from the output of the previous outer iteration, and the quality of this guess improves as the number of iterations increases.}
}
\end{figure}

For Algorithm 2 the initial guess should depend on the iteration $k$. For the first iteration $k=1$ we recommend $\bw_{\text{init}}=\mathbf{0}$. For iterations $k\ge 2$ we recommend $\bw_{\text{init}}=\bw^{(k-1)}$, the solution of the optimisation step from the previous iteration. As the number of iterations increases and the points $\bx_i^{(k)}$ converge, the initial guess $\bw_{\text{init}}=\bw^{(k-1)}$ becomes better and better and consequently the optimisation step becomes quicker and quicker. \change{This is illustrated in Figure \ref{Fig:RelRunTimes}, which shows the relative run time of each iteration for an example in which the relative sizes and relative proportions of small and large grains are the same as in Example \ref{Ex:2} but the number of grains is $n=500$.}
We see that the total runtime of the algorithm is not proportional to the number of iterations $K$; most of the expense is in the first few iterations.

Note that for the first iteration, $k=1$, the initial guess $\bw_{\text{init}}=\mathbf{0}$ does not incorporate any information about the locations $\bx_i^{(0)}$.
It is possible to improve the speed of the first iteration by using a more sophisticated choice of $\bw_{\text{init}}$, e.g., using the \emph{multilevel methods} of M\'erigot \cite{Merigot2011} and L\'evy  \cite{LevySemiDiscrete2015},
which generate a good initial guess $\bw_{\text{init}}$ by solving a sequence of smaller optimisation problems with fewer grains. (For example, you can obtain a good initial guess $\bw_{\text{init}}$ for $n$ grains by first solving a coarser problem with $n/2$ grains;
you can obtain a good initial guess for $n/2$ grains by first solving a coarser problem with $n/4$ grains, etc.)
We found that M\'erigot's multilevel method \cite{Merigot2011} in 2D could halve the run time of iteration $k=1$ when there are $n=10,000$ grains.
It is reported that L\'evy's multilevel program GEOGRAM can handle one million grains in 3D \cite[Table 4]{LevySemiDiscrete2015}.

It is also possible to obtain a better initial guess $\bw_{\text{init}}$ for iterations $k \ge 2$ as follows. The Lloyd step \eqref{eq:RegStep} of Algorithm 2 could be replaced with a \emph{damped Lloyd step} of the form
\begin{equation}
\bx_i^{(k)} := (1-\lambda) \bx_i^{(k-1)}+\lambda \, \frac{1}{\big|L_i^{(k-1)}\big|} \int_{L_i^{(k-1)}} \bx \, \d \bx
\end{equation}
where $\lambda$ is a damping parameter between $0$ and $1$. The choice $\lambda=1$ corresponds to the Lloyd step \eqref{eq:RegStep}. The closer $\lambda$ is to $0$, the closer $\bx_i^{(k)}$ is to $\bx_i^{(k-1)}$, and so the better the \change{associated} initial guess $\bw_{\text{init}}=\bw^{(k-1)}$. Therefore the optimisation step is faster for smaller $\lambda$. On the other hand, the regularisation step has less effect for smaller $\lambda$, and it is necessary to increase the number of iterations $K$ to achieve the same amount of regularisation. For our purposes the full Lloyd step $\lambda=1$ was sufficiently fast and so we did not try to optimise the choice of $\lambda$.

\subsubsection{Choice of the tolerance}
For simplicity we chose the tolerance $\varepsilon$ of the optimisation step of Algorithm 2 to be fixed at each iteration $k$ (recall that the optimisation step terminates when \change{$\big|\nabla g_k\big(\bw^{(k)}\big)\big| < \varepsilon \min_j m_j$)}. The algorithm could be \rev{sped} up, however, by taking $\varepsilon=\varepsilon_k$ to depend on $k$. In order for Algorithm 2 to produce a Laguerre diagram with grains of given volumes up to
\change{a relative error of $\varepsilon$}, we only need the tolerance to be $\varepsilon$ \emph{at the final iteration}, $\varepsilon_K=\varepsilon$. For previous iterations we could use a cruder tolerance: $\varepsilon=\varepsilon_K < \varepsilon_{K-1} < \cdots < \varepsilon_2 < \varepsilon_1$. It is tempting to think that the larger the tolerance, the faster the optimisation step. On the other hand, if $\varepsilon_{k-1}$ is \change{larger} than $\varepsilon_k$, then the initial guess $\bw_{\text{init}}=\bw^{(k-1)}$ at iteration $k$ may be worse, and the optimisation step at iteration $k$ may be slower. So the tolerances $\varepsilon_k$ must be chosen carefully. The choice of fixed tolerance $\varepsilon_k=\varepsilon$ for all $k$ is a simple, reliable option, which is why we used it.

\subsubsection{Choice of the optimisation algorithm}
The speed of the optimisation step depends of course not only on the choice of the initial guess $\bw_{\text{init}}$ and the tolerance $\varepsilon$, but also on the choice of the optimisation algorithm. For example, instead of using a quasi-Newton method like we did, one could use
Newton's method. Newton's method tends to converge faster than quasi-Newton methods (quadratically rather than superlinearly), although it is harder to implement since it requires the second derivative of $g$ (whereas quasi-Newton methods only require the first derivative) \cite{BoydVandenberghe}.

It can be shown (see, e.g., \cite{BourneRoper15}) that
\begin{equation}
\frac{\partial^2 g}{\partial w_i \partial w_j} = -\frac{\partial |L_i|}{\partial w_j} =\begin{cases}
\displaystyle
-\sum_{k\in N_i} \frac{a_{ik}}{2|\bx_i-\bx_k|} & \text{if $i=j$},\\[1.6em]
\displaystyle
\frac{a_{ij}}{2|\bx_i-\bx_j|} & \text{if $j\in N_i$},\\[1.6em]
\displaystyle
0 & \text{otherwise},
\end{cases}
\end{equation}
where $a_{ij}$ is the area of the face between cell $i$ and cell $j$ and $N_i$ is the index set of the neighbours of cell $i$ (that is $j\in N_i$ if and only if cell $j$ and cell $i$ share a face). \change{The pseudo-inverse of this Hessian matrix is used in a damped Newton method in \cite[Algorithm 1]{KitagawaMerigotThibert2019},
where fractions of a full Newton step are used
in order to control the error reduction and the minimum volume of a cell (to stop cells disappearing). The authors prove that their damped Newton method converges
globally with order $1$ and locally with order $2$ \cite[Theorem 1.5]{KitagawaMerigotThibert2019}.}

\subsection{Initialisation: Effect on the spatial distribution }
In this section we discuss the initialisation step of Algorithms 1 and 2.

\subsubsection{Initialisation of the seeds}
\label{Subsec:InitSeeds}
\begin{figure}[tbp]
  \begin{subfigure}[b]{0.5\textwidth}
    \centering
    \subcaption{\label{sfig:random}Random.}
    \includegraphics[trim={0cm 0.25cm 0cm 0.5cm},clip,width=\textwidth]{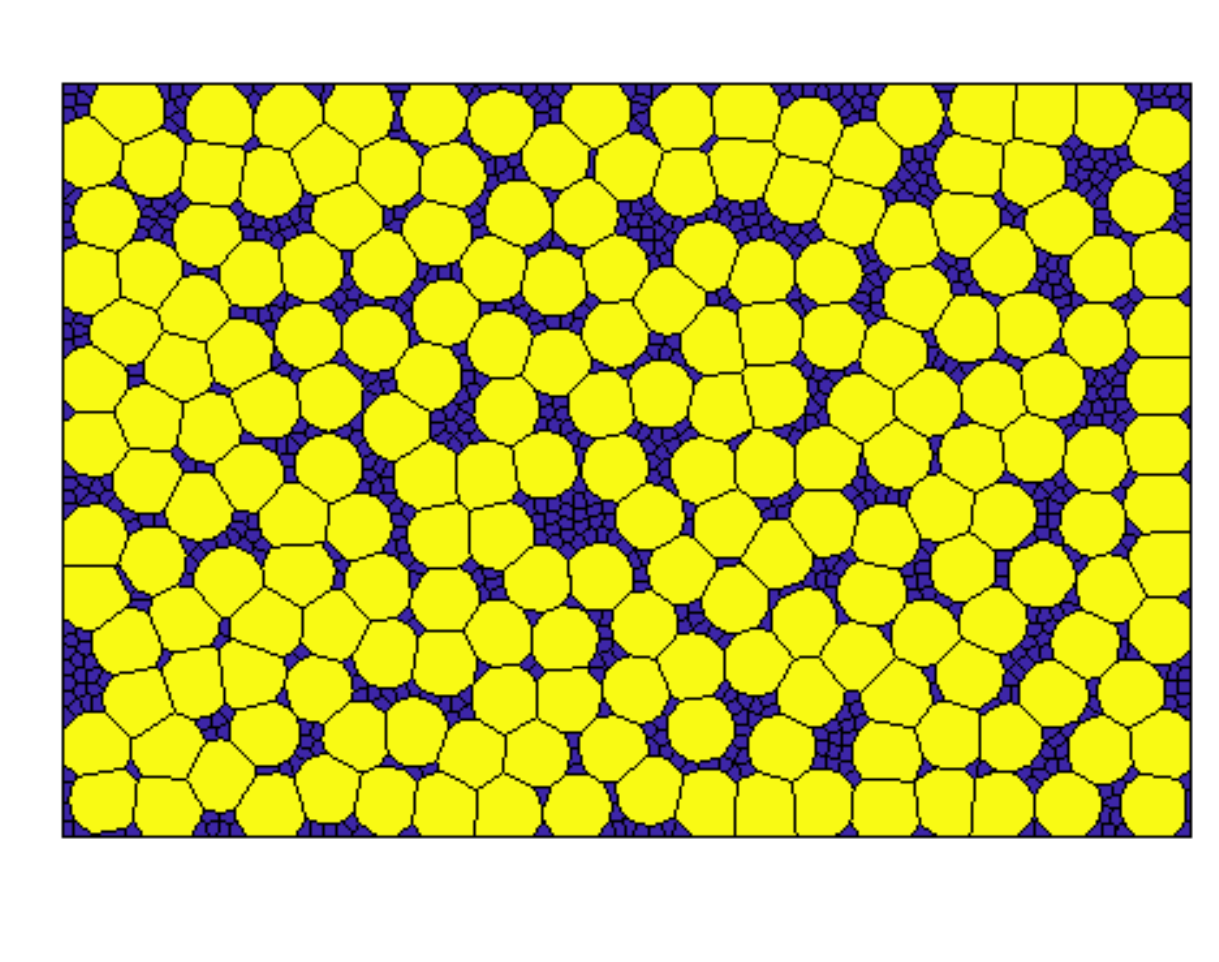}
  \end{subfigure}
  \begin{subfigure}[b]{0.5\textwidth}
    \centering
    \subcaption{\label{sfig:banded}Banded.}
    \includegraphics[trim={0cm 0.25cm 0cm 0.5cm},clip,width=\textwidth]{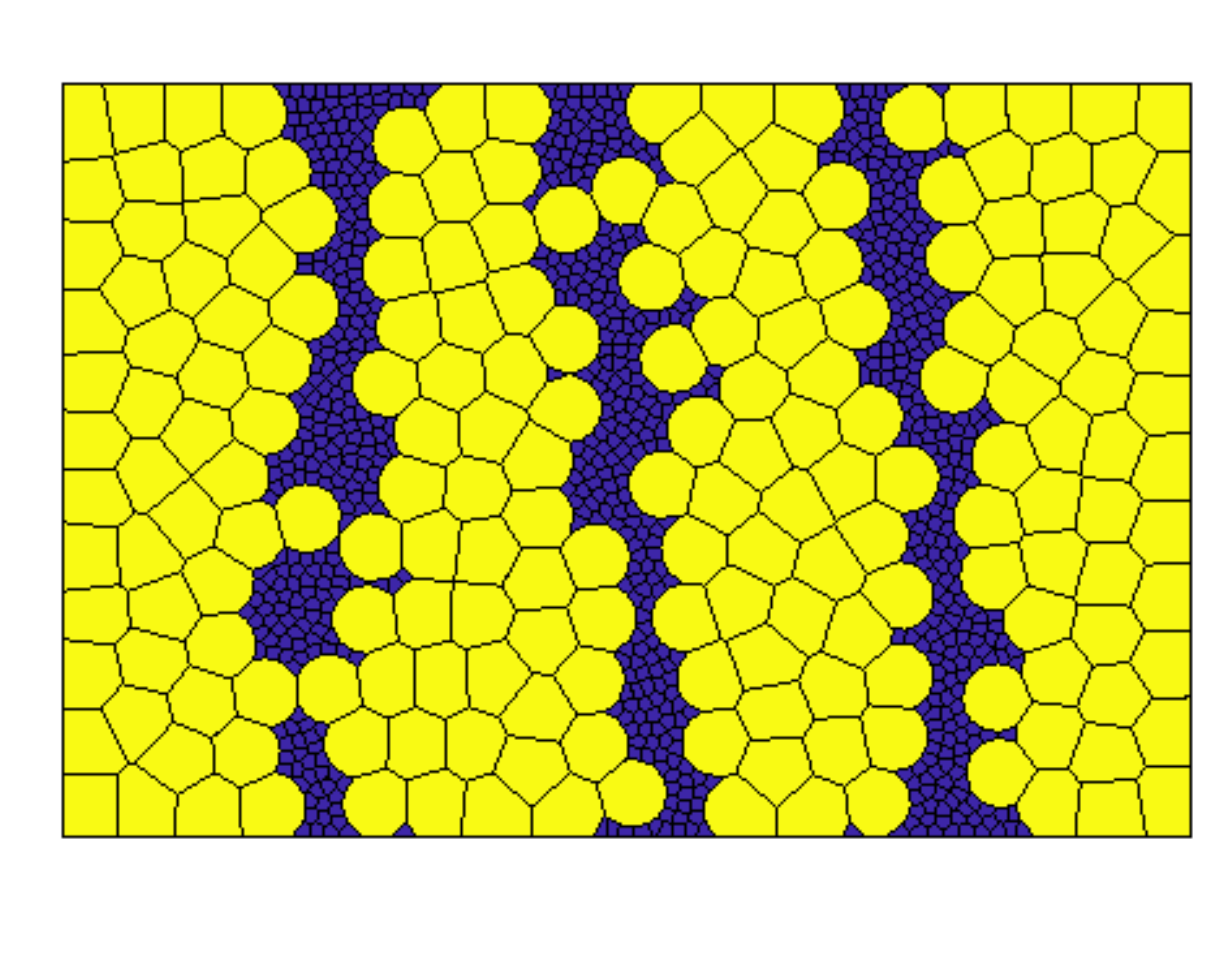}
  \end{subfigure}

  \begin{subfigure}[b]{0.5\textwidth}
    \centering
    \subcaption{\label{sfig:cluster}Clustered.}
    \includegraphics[trim={0cm 0.25cm 0cm 0.5cm},clip,width=\textwidth]{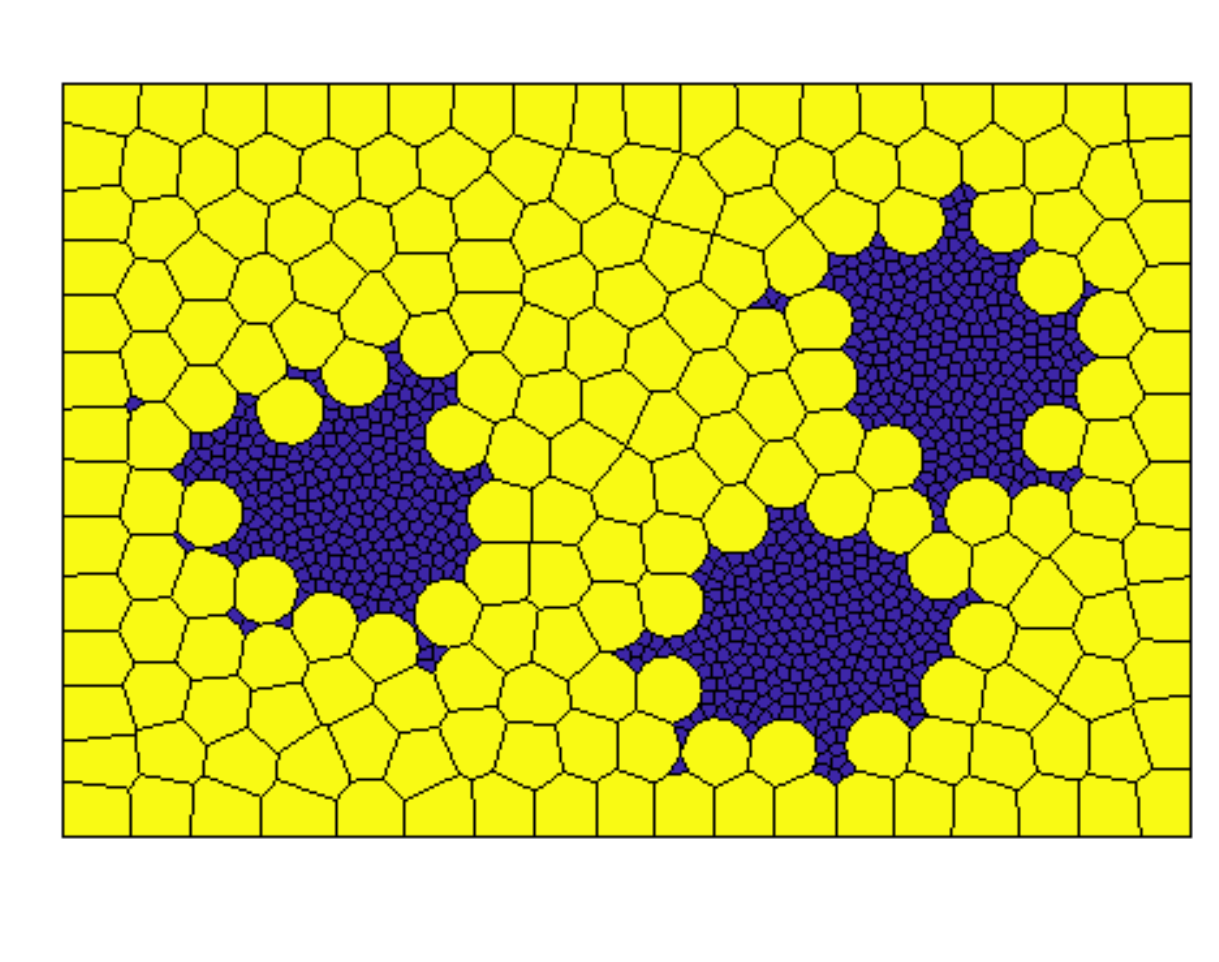}
  \end{subfigure}
  \begin{subfigure}[b]{0.5\textwidth}
    \centering
    \subcaption{\label{sfig:mixed}Mixed banded and random.}
    \includegraphics[trim={0cm 0.25cm 0cm 0.5cm},clip,width=\textwidth]{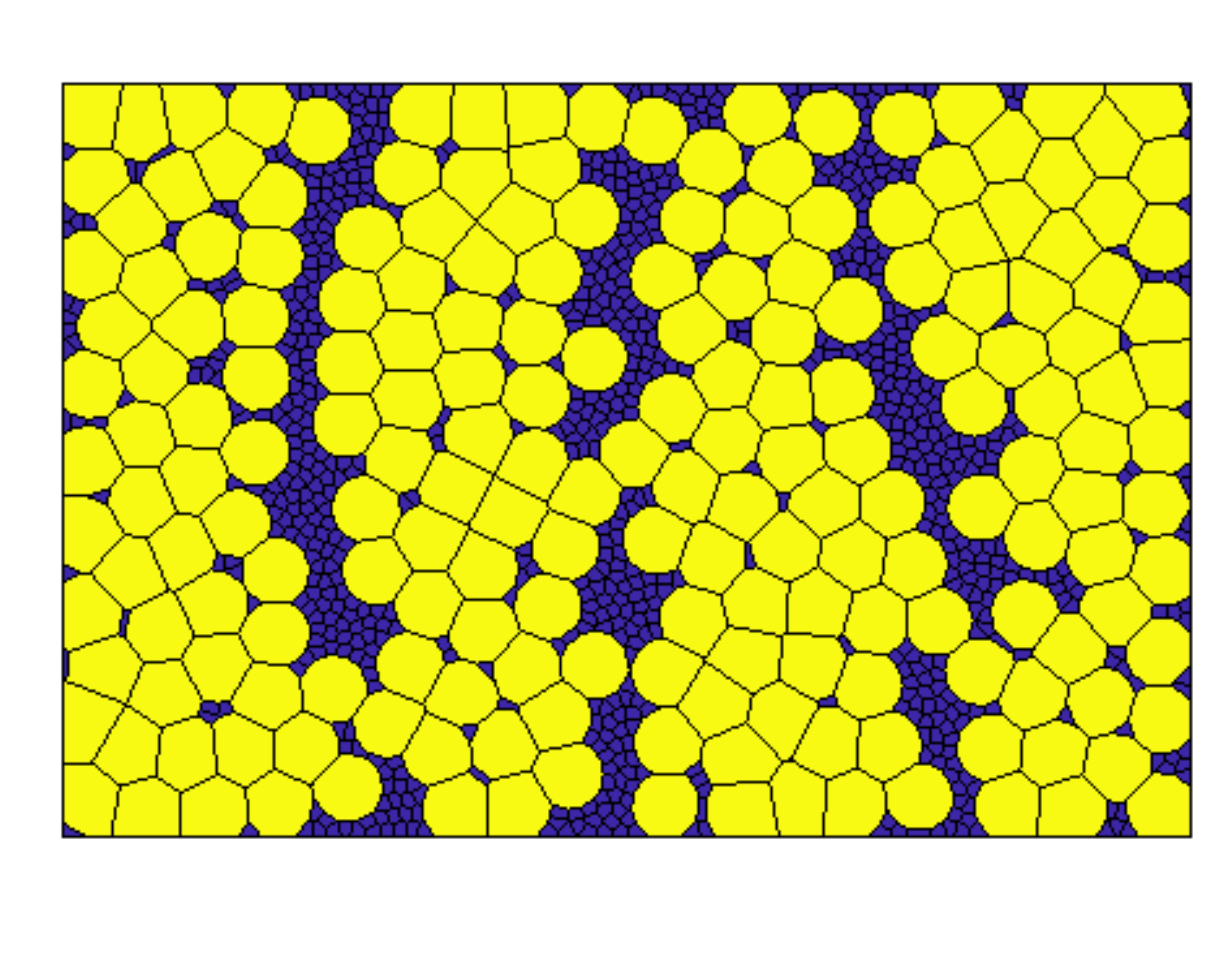}
  \end{subfigure}

  \caption{\label{fig:init_seeds_bimodal}
  \todo{Examples of more advanced microstructures: Coupling of size and spatial distributions with Algorithm 2.
  This figure shows} the output of Algorithm 2 after $K=20$ iterations with different initial distributions for the seed locations. In all cases there are $n=1000$ grains with $n_1=800$ grains of size $x$ (in dark blue) and $n_2=200$ grains of size $20x$ (in yellow).  \subref{sfig:random} Random distribution of initial locations. Here the initial generator locations of the large and small grains are uniformly distributed over $\Omega$. \subref{sfig:banded} Banded distribution of initial locations. Here the different sized grains have initial generator locations that lie inside bands within $\Omega$. The sizes of the bands have been chosen so that there are approximately equal numbers of small grains within each blue band and approximately equal numbers of large grains within each yellow band. \subref{sfig:cluster} Clustered distribution of initial locations. Here the smaller grains have initial generator locations that lie inside non-overlapping discs. \subref{sfig:mixed} A mixed distribution: the initial generators are such that the larger grains are arranged in bands and the smaller grains are a combination of the banded and random distributions.}
\end{figure}

The locations of the generators $\bx_1,\ldots,\bx_n$ at the termination of Algorithm 2 is a strong function of the initial choice $\bx_1^{(0)},\ldots,\bx_n^{(0)}$. This simple observation gives
us a great deal of control over the spatial distribution of the different sized grains.
 Examples of Algorithm 2 with different initial distributions of the generators $\bx_1^{(0)},\ldots,\bx_n^{(0)}$ are shown in Figure \ref{fig:init_seeds_bimodal}. In these examples $\Omega=[0,3]\times[0,2]$ and there are $n=1000$ grains. There are $n_1=800$ grains of size $x$ and $n_2=200$ grains of size $20x$. The tolerance is \change{$\varepsilon=0.01$} and the number of iterations of Algorithm 2 is $K=20$. The figure shows the output
of Algorithm 2. The final spatial distribution of grains has some features in common with the spatial distribution of the initial generator locations.

\begin{figure}[!htbp]
  \begin{subfigure}[t]{0.5\textwidth}
    \centering
    \subcaption{\label{sfig:gradient} Increasing gradient of grain size.}
    \includegraphics[trim={0cm 0.25cm 0cm 0.5cm},clip,width=\textwidth]{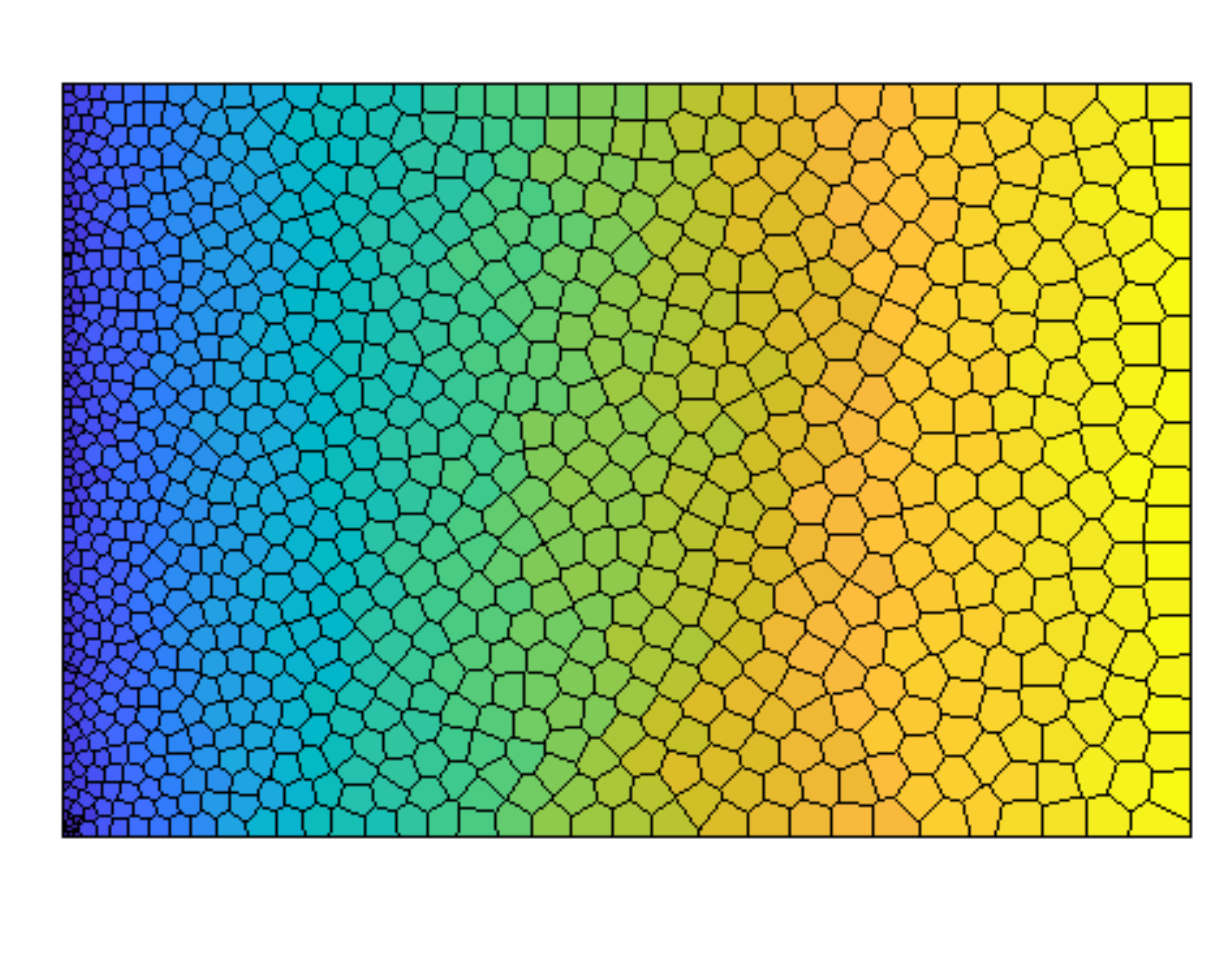}
  \end{subfigure}
  \begin{subfigure}[t]{0.5\textwidth}
    \centering
    \subcaption{\label{sfig:doublegradient} Varying gradient of grain size.}
    \includegraphics[trim={0cm 0.25cm 0cm 0.5cm},clip,width=\textwidth]{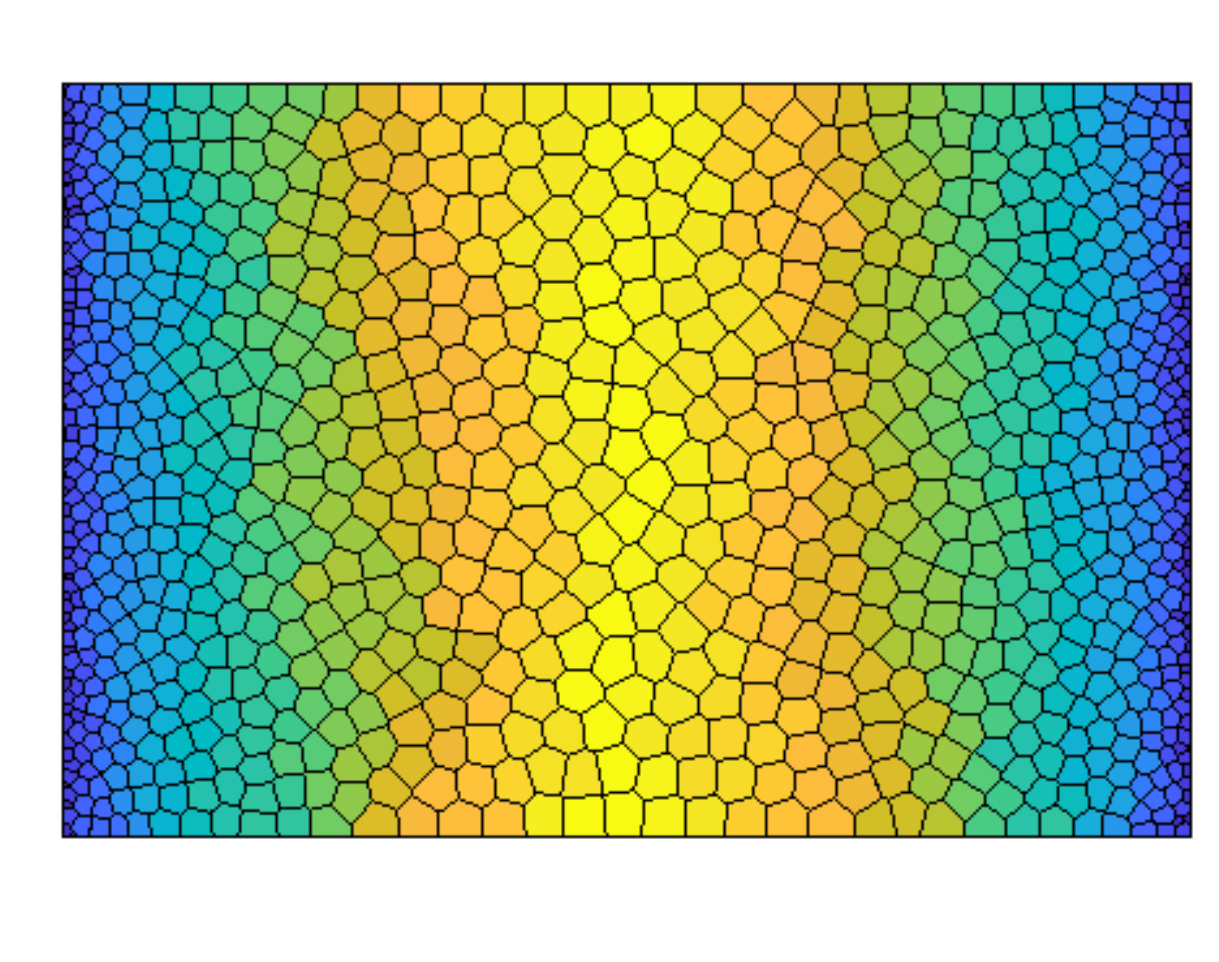}
  \end{subfigure}
\caption{\label{fig:init_seeds_gradient} Experiments to demonstrate a gradient in the distribution of grain sizes. Here $n=1000$ grains have areas drawn from a uniform distribution such that the ratio of the largest to the smallest grain size is at most $100$. The domain is $\Omega=[0,3]\times[0,2]$. \change{The results are shown after $K=20$ iterations and the grains are coloured according to their area.} \subref{sfig:gradient} The initial seed locations $\bx_1^{(0)},\ldots, \bx_n^{(0)}$ are distributed such that the $x$-coordinate increases with grain size. \subref{sfig:doublegradient} The initial seed locations are distributed such that the larger grains are found in the middle of $\Omega$.}
\end{figure}

A further example of controlling the spatial distribution of grains can be seen in Figure \ref{fig:init_seeds_gradient}. In this example $n=1000$ grains have areas drawn from a random distribution such that the ratio of  the largest to the smallest grain size is at most $100$. The Laguerre diagram in Figure \ref{fig:init_seeds_gradient}(a) has the property that the grain sizes tend to increase \change{from} left to right.
A variety of spatial distributions of grain sizes can be simulated by first distributing the seed locations appropriately. Figure \ref{fig:init_seeds_gradient}(b) shows how a more complicated distribution can be produced.

\subsubsection{Initialisation of the weights}
The choice of $\bw^{(0)}$ in the initialisation step of Algorithm 2 is also important. One should choose $w_1^{(0)},\ldots,w_n^{(0)}$ so that the Laguerre diagram generated by $\big(\bx_1^{(0)},w_1^{(0)}\big),\ldots,\big(\bx_n^{(0)},w_n^{(0)}\big)$ has no empty Laguerre cells. If there are empty cells then
the regularisation step is not defined (there is division by zero in equation \eqref{eq:RegStep} if $L_i^{(0)}$ is empty). A good choice is $\bw^{(0)}=\mathbf{0}$ since then the Laguerre diagram generated by $(\bx_1^{(0)},w_1^{(0)}),\ldots,(\bx_n^{(0)},w_n^{(0)})$ is a Voronoi diagram and so has no empty cells, whatever the choice of $\bx_1^{(0)},\ldots,\bx_n^{(0)}$.

\subsection{Stopping criteria}
\label{Subsec:StopCrit}
In Algorithm 2 the user must specify the number of regularisation steps $K$.
As we discussed in Section \ref{Sec:Algorithms}, for large values of $K$ the Laguerre diagram \change{resulting from} Algorithm 2 is a approximately a \emph{centroidal Laguerre diagram}, which means that each seed $\bx_i^{(K)}$ is approximately the centre of mass of its Laguerre cell $L_i^{(K)}$. Centroidal Laguerre diagrams tend to have very regular-shaped cells, e.g., in 2D if the grains all have the same target areas, $m_i=1/n$, and if $n$ and $K$ are large, then the Laguerre diagram looks locally like a regular hexagonal tiling.
Indeed for steel microstructures we found that if $K$ is too large, then Algorithm 2 tends to produce grains that are too `round' compared to EBSD measurements of grain aspect ratios.

Instead of fixing the number of steps $K$ in advance, the user could terminate the algorithm whenever some measure of the maximum grain aspect ratio falls below a critical threshold.
\change{For example, the aspect ratio of a grain can be measured using the ratio of its largest and smallest principal moments of inertia, or using the ratio of the radii of circumscribed and inscribed balls,
	or using its sphericity \cite{Spettl2014}, which is the ratio of the surface area of the volume-equivalent sphere to the surface area of the grain.}

\change{The user may also want} to terminate the algorithm if the distance $| \bx_i^{(k)}-\bx_i^{(k-1)} |$ moved by the seeds from one iteration to the next falls below some threshold.
The Laguerre diagram $\{ L_i^{(k)} \}$ tends to evolve slowly with $k$ when $k$ is large, as illustrated in Figure \ref{Fig:Alg2}, and the evolution can slow down dramatically when there is a T1 topological transition (to borrow a term from foam dynamics). This topological transition involves a change of cell neighbour relations; in 2D this is via \change{coalescence} of two triple junctions of cell boundaries. So in general there is little to be gained from performing a large number of regularisation steps, especially since our aim is not to produce a centroidal Laguerre diagram, but rather to produce a physically realistic microstructure. \change{(If on the other hand the user's aim is to produce centroidal Laguerre diagrams, then it would be better to use the quasi-Newton method of \cite{LevyCentroidalPower}, which converges superlinearly as opposed to the linear convergence of Lloyd's algorithm.)}


\section{Examples}
\label{Sec:Examples}
This section includes some large examples in 3D to illustrate the capabilities of our algorithms.

\begin{example}[Run time tests]
	\label{Example:RunTimeTests}
	\change{Figure \ref{Fig:RunTimes} gives some run times of Algorithm 2 for creating monodisperse and polydisperse periodic Laguerre diagrams in 3D.
		Here Algorithm 2 has been used to create $n/2$ grains of volume $x$ and $n/2$ grains of volume $rx$, for $r=1$ and $r=5$, in a cube of side length $100$ with error tolerance $\varepsilon=0.01$.
		We see that the run time grows roughly quadratically in the range $n=2000$ to  $n=5000$,
		for both the monodisperse (single phase) case $r=1$ and the polydisperse (dual phase) case $r=5$. 
		This could be expected since the cost of each iteration of the BFGS method used by \emph{fminunc} is $\mathcal{O}(n^2)$ (a matrix-vector multiplication). Also the worst-case optimal time it takes to compute a Laguerre diagram of $n$ cells is $\mathcal{O}(n^2)$ in 3D, although we found that in these examples the cost of computing the Laguerre diagrams was sub-quadratic (but super-linear); see also the discussion in Section \ref{Subsec:Lag}. For $n>5000$ the run time grows a little faster than $\mathcal{O}(n^2)$.
		Observe also from Figure \ref{Fig:RunTimes} that it is \rev{about 50\%} more expensive to compute dual phase RVEs ($r=5$) than single phase RVEs ($r=1$).}
	
	\begin{figure}[!hbtp]
		\centering
		\begin{tabular}{ccc}
			\hline
			$n$ & \multicolumn{2}{c}{run time (s)}
			\\
			\hline
			& $r=1$ & $r=5$
			\\
			1000 & 1.96 & 2.86 \\
			2000 & 5.26 & 8.41 \\
			3000 & 11.41 & 17.62 \\
			5000 & 31.53 & 46.84 \\
			7500 & 75.86 & 110.09 \\
			10000 & 147.15 & 206.71 \\
			15000 & 356.41 & 525.33 \\
			20000 & 689.98 & 989.42 \\
			\hline
		\end{tabular}
		\quad
		\includegraphics[align=c,width=0.6\textwidth]{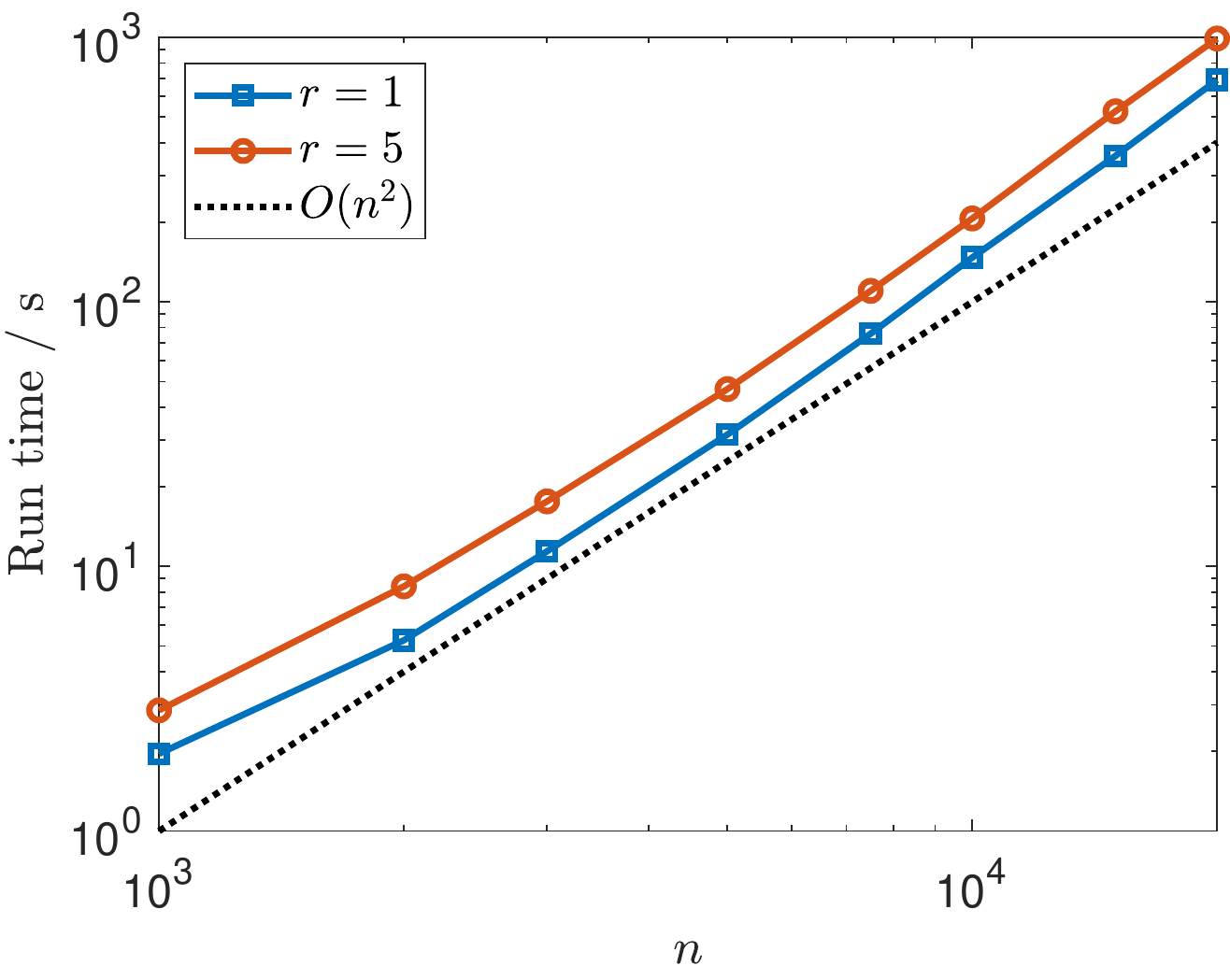}
		\caption{\label{Fig:RunTimes} Run times in seconds of Algorithm 2 for creating monodisperse (single phase) and polydisperse (dual phase) periodic RVEs in 3D. There are $n/2$ grains of volume $x$ and $n/2$ grains of volume $rx$ in a cube of side length $100$ with at most \change{$1\%$ error ($\varepsilon=0.01$)} in the volumes of the grains ($r=1$ corresponds to a single phase material, $r=5$ corresponds to a dual phase material, $x$ is chosen \change{such} that the total volume of the grains equals the volume of the box). We used $K=5$ regularisation steps, and the initial seed locations $\bx_1^{(0)},\ldots,\bx_n^{(0)}$ were chosen randomly from a uniform distribution. The simulations were performed on an Intel Xeon E5-1620V3 (3.5GHz, 4 cores, 8 threads). The graph on the right has a log-log scale. The black dotted line is the graph of the function $cn^2$, where $c$ is a constant.
			It is included to illustrate how the run times grow with $n$.}
	\end{figure}
\end{example}

\begin{figure}[!htbp]
  \centering
  \includegraphics[align=c,width=\textwidth]{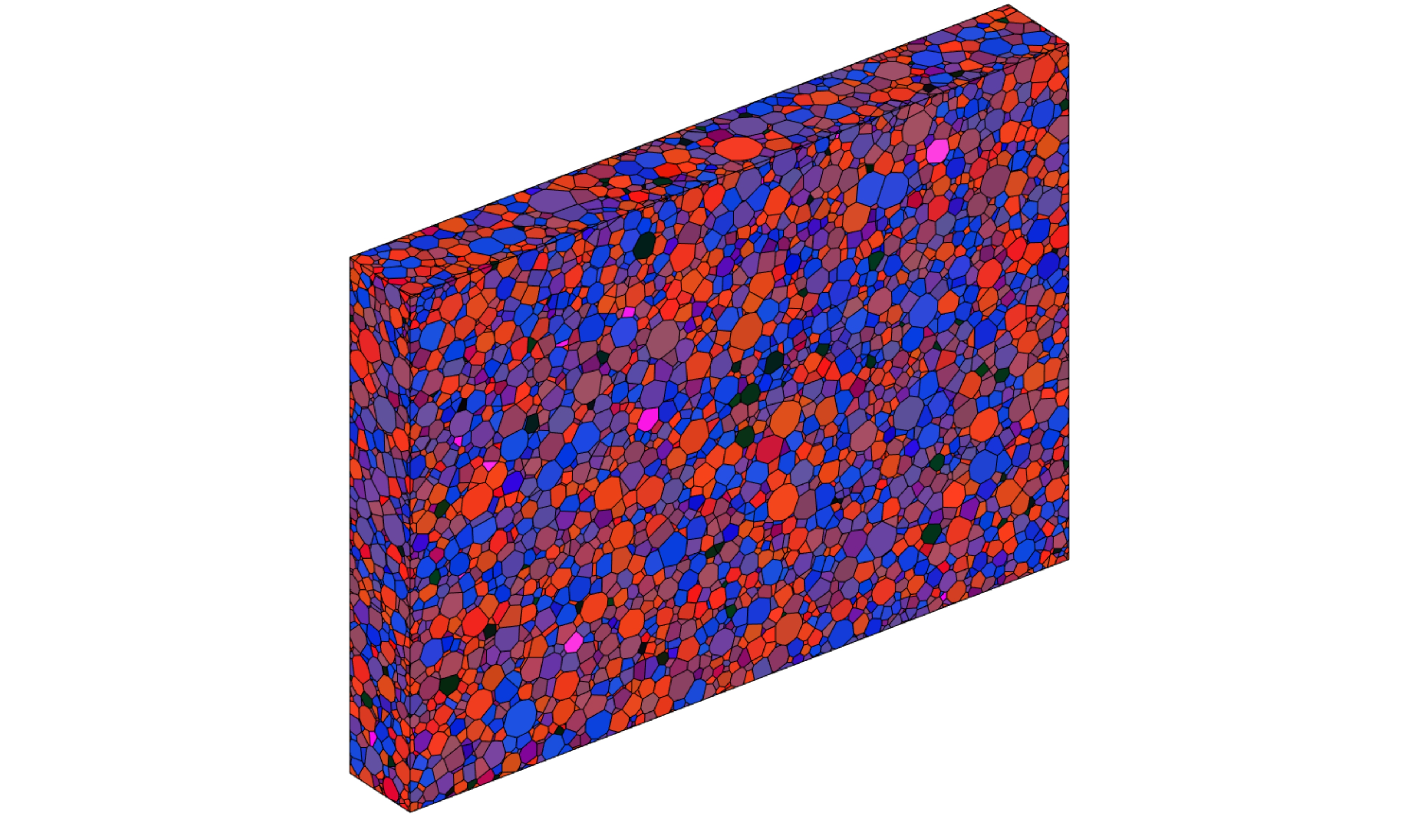}
  \caption{\label{fig:EBSDperiodic} A periodic RVE of an IF steel with $n=9211$ grains, generated using Algorithm 2 (see Example \ref{Example:IFsteelPeriodicRVE}). The grains are coloured according to their lattice \rev{orientations}.}
\end{figure}

\begin{figure}[!htbp]
  \centering
 \includegraphics[width=0.6\textwidth]{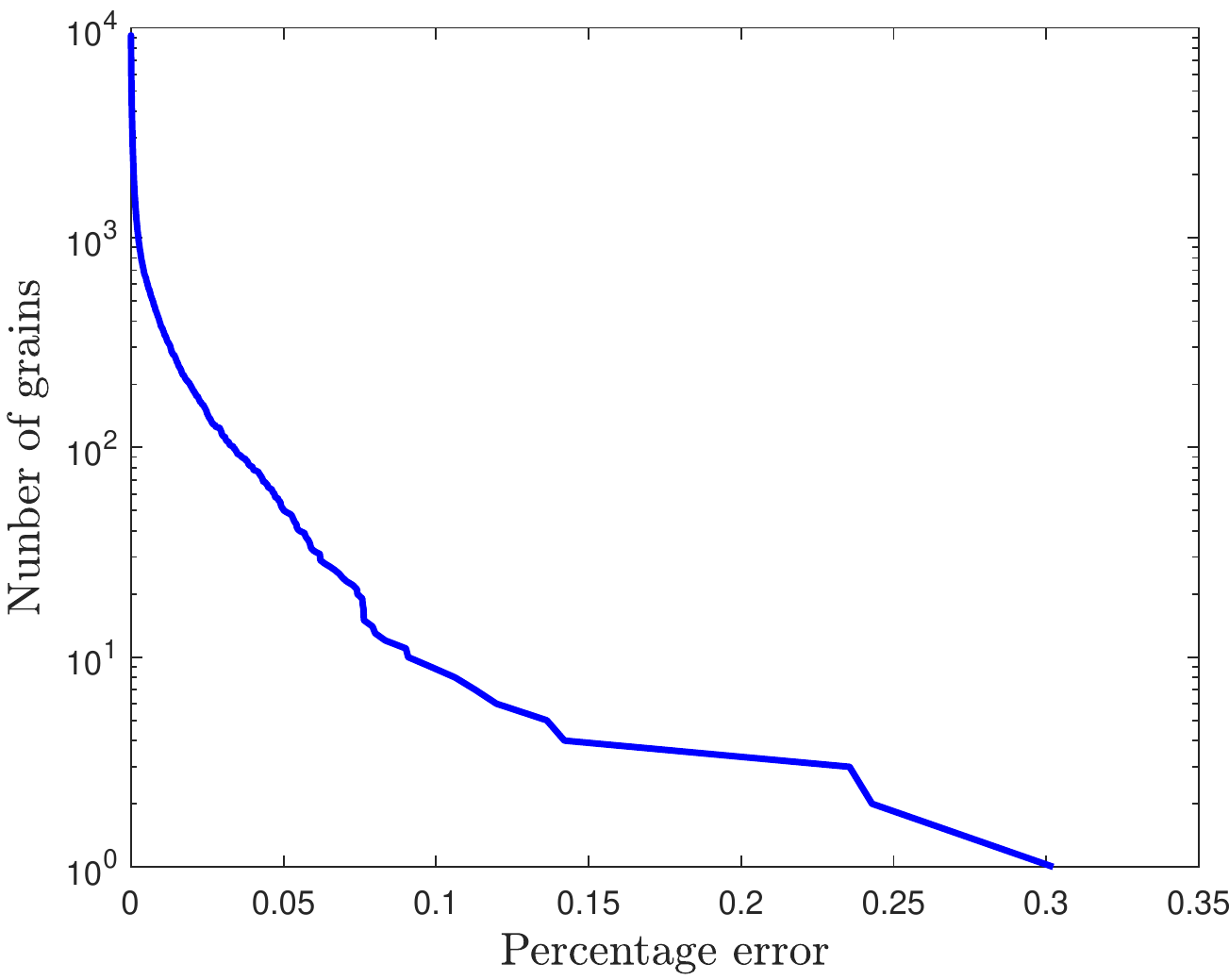}
  \caption{\label{fig:EBSDperiodic_errors}
  	\change{The complementary cumulative number distribution of the percentage error of the volumes of the grains for Example \ref{Example:IFsteelPeriodicRVE}. For a percentage error $x$ we plot the number of grains that have volume percentage error at least $x$. Most of the grain volumes have percentage error less than $0.05\%$, which is an order of magnitude below the tolerance of \change{$0.5\%$ ($\varepsilon=0.005$)}. In this example, the maximum percentage error is $0.30\%$.}}
\end{figure}

\begin{example}[Generating a periodic RVE of an IF steel]
\label{Example:IFsteelPeriodicRVE}
Figure \ref{fig:EBSDperiodic} shows an example of a periodic Laguerre diagram created using Algorithm 2. The target volumes $m_i$ are taken from a 3D EBSD measurement of an IF (interstitial free) steel (EN 10130 grade DC06). There are $n=9211$ grains in a box of dimensions $670\mu m\times80 \mu m\times480 \mu m$. We took the initial seed locations $\bx_1^{(0)},\ldots,\bx_n^{(0)}$ to be the centres of mass of the grains from the EBSD data, and performed $K=10$ regularisation steps with a tolerance of \change{$0.5\%$ ($\varepsilon=0.005$)}. The grains in Figure \ref{fig:EBSDperiodic} are coloured according to their lattice \rev{orientations} by mapping the three Euler angles linearly to RGB values. The orientations were taken from the EBSD data. Figure \ref{fig:EBSDperiodic_errors} shows that the volumes of all the grains are correct to within $0.5\%$, and most volumes are correct to within $0.05\%$.
\end{example}

\begin{example}[Fitting a Laguerre diagram to EBSD measurements]
\label{Ex:DataFitting}
 The main aim of this paper is to create Laguerre diagrams with given volume distributions for use in computational homogenisation simulations. We briefly mention, however, how Algorithm 1 can be used to fit a Laguerre diagram to EBSD data of grain volumes and centroids. Figure \ref{Fig:Fitting} is an example of a non-periodic Laguerre diagram fitted to a 3D EBSD measurement of an IF steel (EN 10130 grade DC06)
 with $n=9211$ grains in a box of dimensions $670\mu m\times80 \mu m\times480 \mu m$ (this is the same EBSD data used in the previous example). In the initialisation step of Algorithm 1 we took $\bx_1,\ldots,\bx_n$  to be the centroids of the grains from
 the EBSD data.  The target volumes $m_i$ were also taken from the EBSD data. We used a tolerance of \change{$1\%$ ($\varepsilon=0.01$)}. As in the previous example, \change{the grains in Figure \ref{Fig:Fitting} are coloured according to their lattice orientation.}
Observe that Figure \ref{Fig:Fitting} is less regular than Figure \ref{fig:EBSDperiodic}, \change{which is because Algorithm 1 is missing the regularisation steps of Algorithm 2}.

Figure \ref{Fig:Fitting_VolErrors} shows that the volumes of all the grains are correct to within $0.56\%$, and most volumes are correct to within $0.1\%$.

\change{Figure \ref{Fig:CentroidErrors} shows \change{the complementary cumulative number distribution of the relative errors of the centroids.}
  The relative error for grain $i$ is defined by
  \begin{equation}
  \frac{| \bx_i - \bc_i|}{r_i}
\end{equation}
   where $\bx_i$ is the centroid of grain $i$ from the EBSD data, $\bc_i$ is the centroid of the Laguerre cell $L_i$, and $r_i$ is the radius of a sphere of volume $m_i$, where $m_i$ is the target volume of grain $i$. This definition of relative error was proposed by \cite{PetrichEtAl2019}.
  As expected, the centroid errors are higher than the volume errors since Algorithm 1 does not directly try to fit the centroids. To be precise, 
  the optimisation step of Algorithm 1 only minimises the fitting error of the volumes; the centroid fitting error is not minimised (the centroids do not appear in the objective function $g$). However, the centroid measurements are used in the initialisation step of Algorithm 1.} The relative error of 79\% of the grain centroids is less than 1 and the relative error of 93\% of the grain centroids is less than 2.

The run time for this example was 376 seconds on an Intel Xeon E5-1620V3 (3.5GHz, 4 cores, 8 threads) with the initial guess
 $(\bw_{\text{init}})_i=(3m_i/(4\pi))^{2/3}$, which \change{was} inspired by sphere-packing methods \cite{DepriesterKubler2019,LLLFP2011,PFCRMMO2019,WuZhouWangYang2010}.

\begin{figure}[!htbp]
\centering
\includegraphics[width=\textwidth]{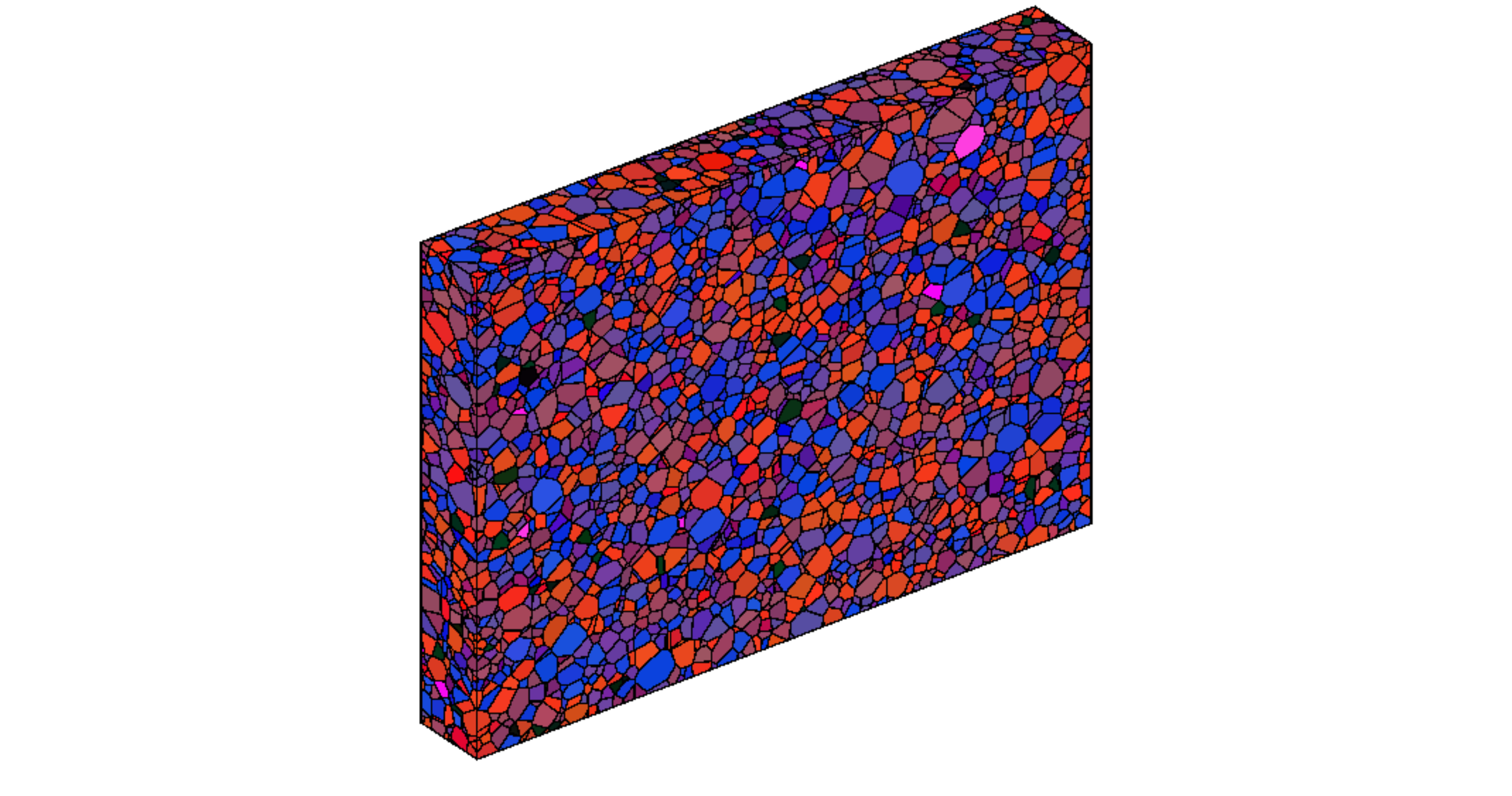}
\caption{\label{Fig:Fitting} A non-periodic Laguerre diagram fitted to 3D EBSD data of an IF steel using Algorithm 1 (see Example \ref{Ex:DataFitting}). The grains are coloured according to their lattice \rev{orientations}. The volume distribution has a fitting error of less than $1\%$. The texture intensity inherits the same fitting error.}
\end{figure}

\begin{figure}[!htbp]
\centering
 \includegraphics[width=0.6\textwidth]{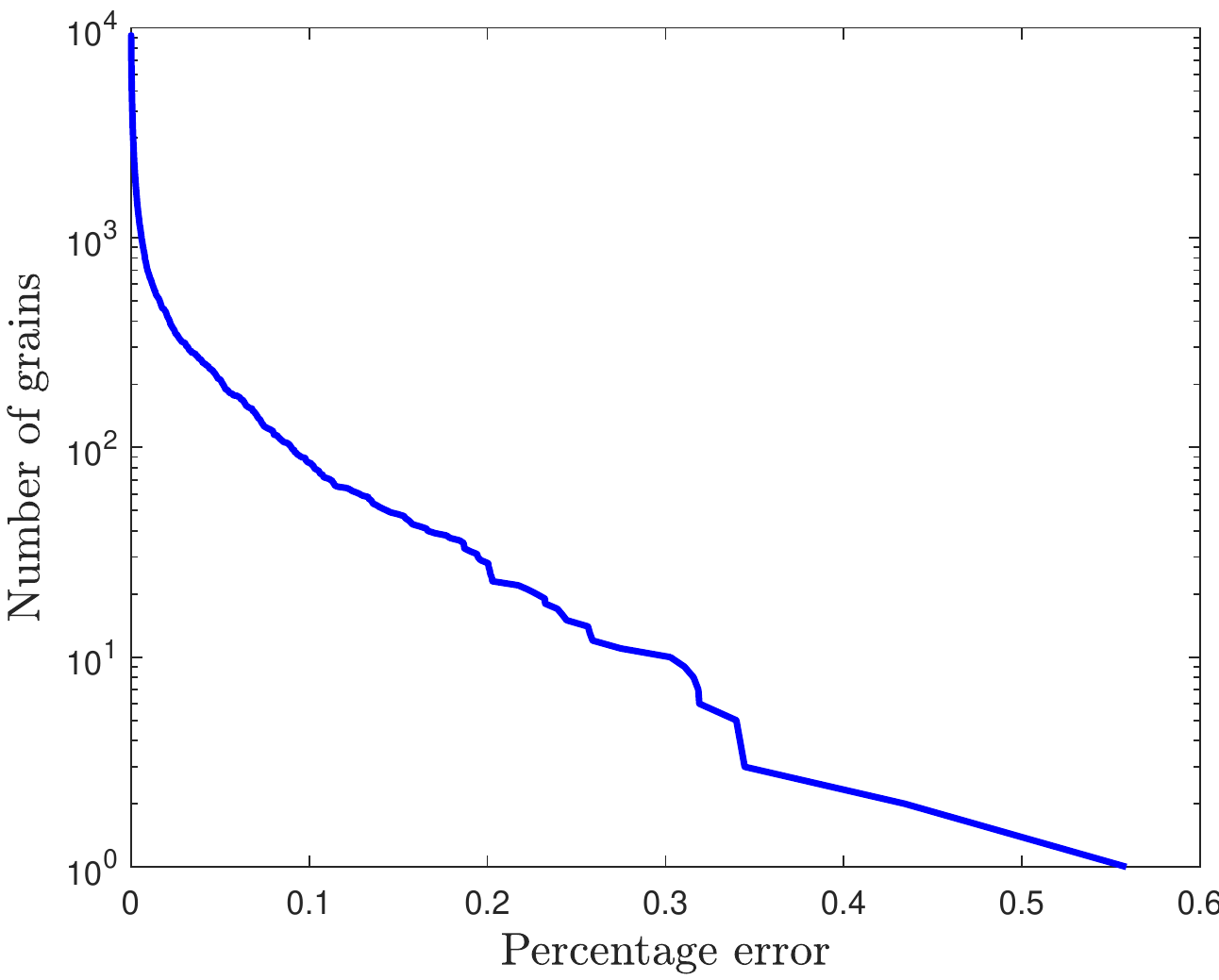}
\caption{\label{Fig:Fitting_VolErrors}
	\change{The complementary cumulative number distribution of the percentage error of the volumes of the grains for Example \ref{Ex:DataFitting}. For a percentage error $x$ we plot the number of grains with volume percentage error at least $x$. The largest percentage error is $0.56\%$ and the second largest is $0.43\%$. All the other percentage errors are below $0.34\%$.}}
\end{figure}

\begin{figure}[!htbp]
\centering
\includegraphics[width=0.6\textwidth]{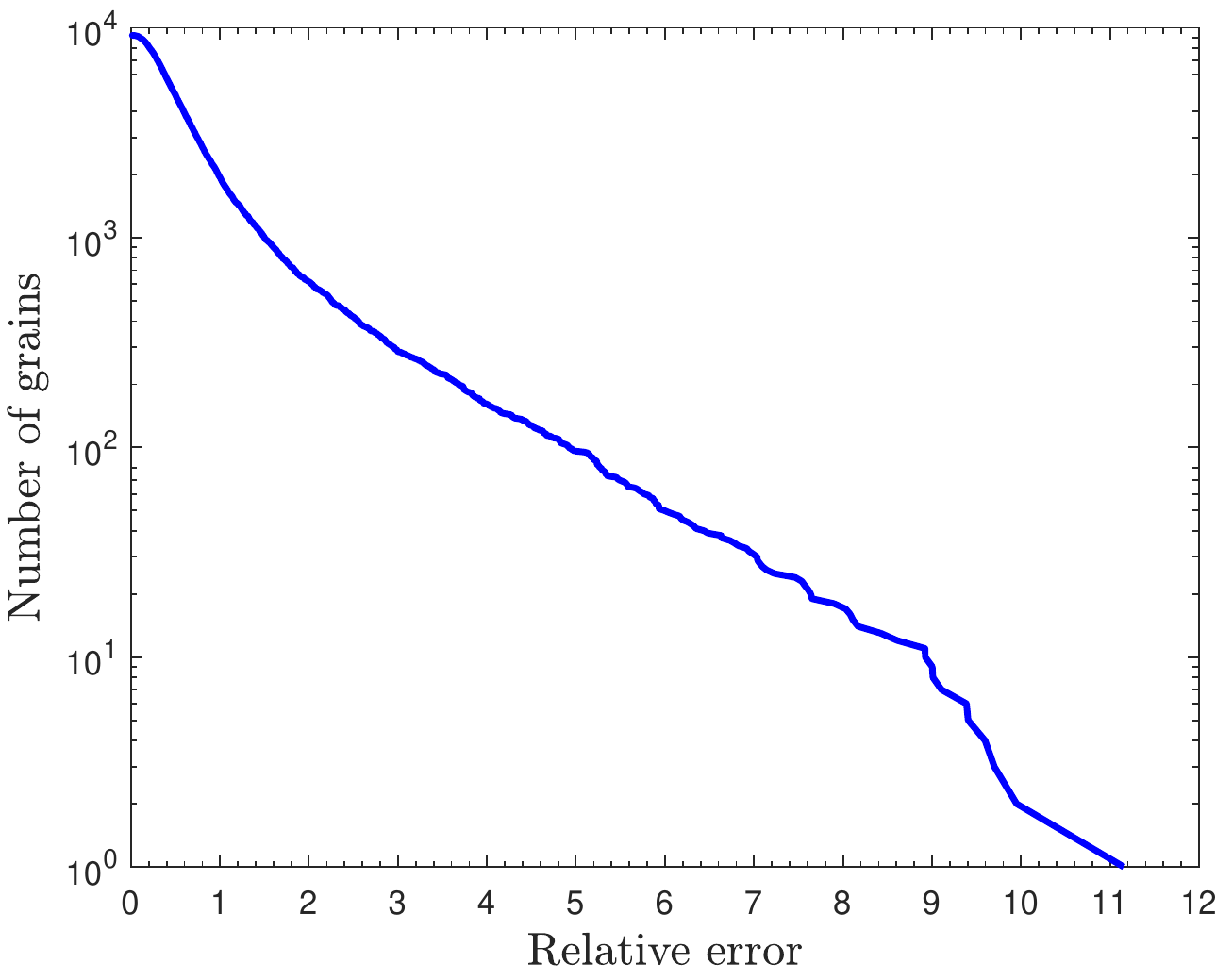}
\caption{\label{Fig:CentroidErrors} \change{The complementary cumulative number distribution of the relative error of the centroids of the grains for Example \ref{Ex:DataFitting}. For a relative error $x$ we plot the number of grains with centroid relative error at least $x$. The maximum relative error is 11.16\%, which is worse than the result given in \cite[Fig.~9]{PetrichEtAl2019}, although for most of the grains the relative errors are comparable: 7278 of the 9211 grains have relative error less than 1, and 8596 of the 9211 grains have relative error less than 2.}}
\end{figure}
\end{example}

\begin{example}[Generating a dual phase RVE with a banded microstructure]
\label{Example:Bands}
Figure \ref{Fig:Bands} shows an example of a periodic, dual phase Laguerre diagram with a band of small grains in the centre, generated using Algorithm 2. There are $n=10,000$ grains: $8000$ grains of volume $x$ and $2000$ grains of volume $20x$ (where $x$ was chosen so that the total volume of the grains equals the volume of the box). We used $K=20$ regularisation steps and a volume tolerance of $1\%$. The grains are coloured according to their volume.
In order to obtain the banded structure, we placed the initial seeds $\bx_1^{(0)},\ldots,\bx_n^{(0)}$ at random within bands of the correct volume. We see from Figure \ref{Fig:Bands} that these bands were largely preserved by the regularisation steps.

\begin{figure}[!htbp]
\centering
\includegraphics[width=0.49\textwidth]{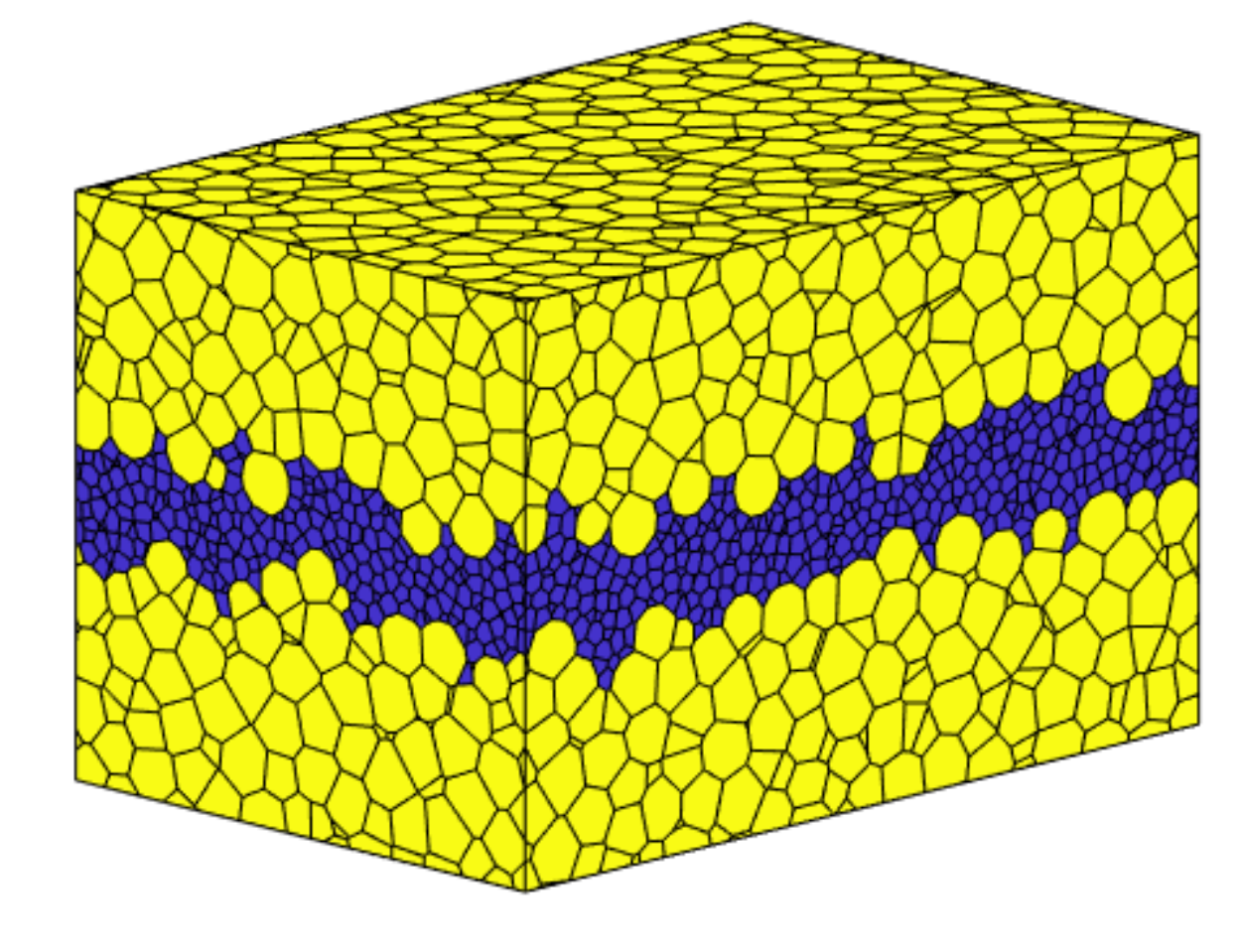}
\includegraphics[width=0.49\textwidth]{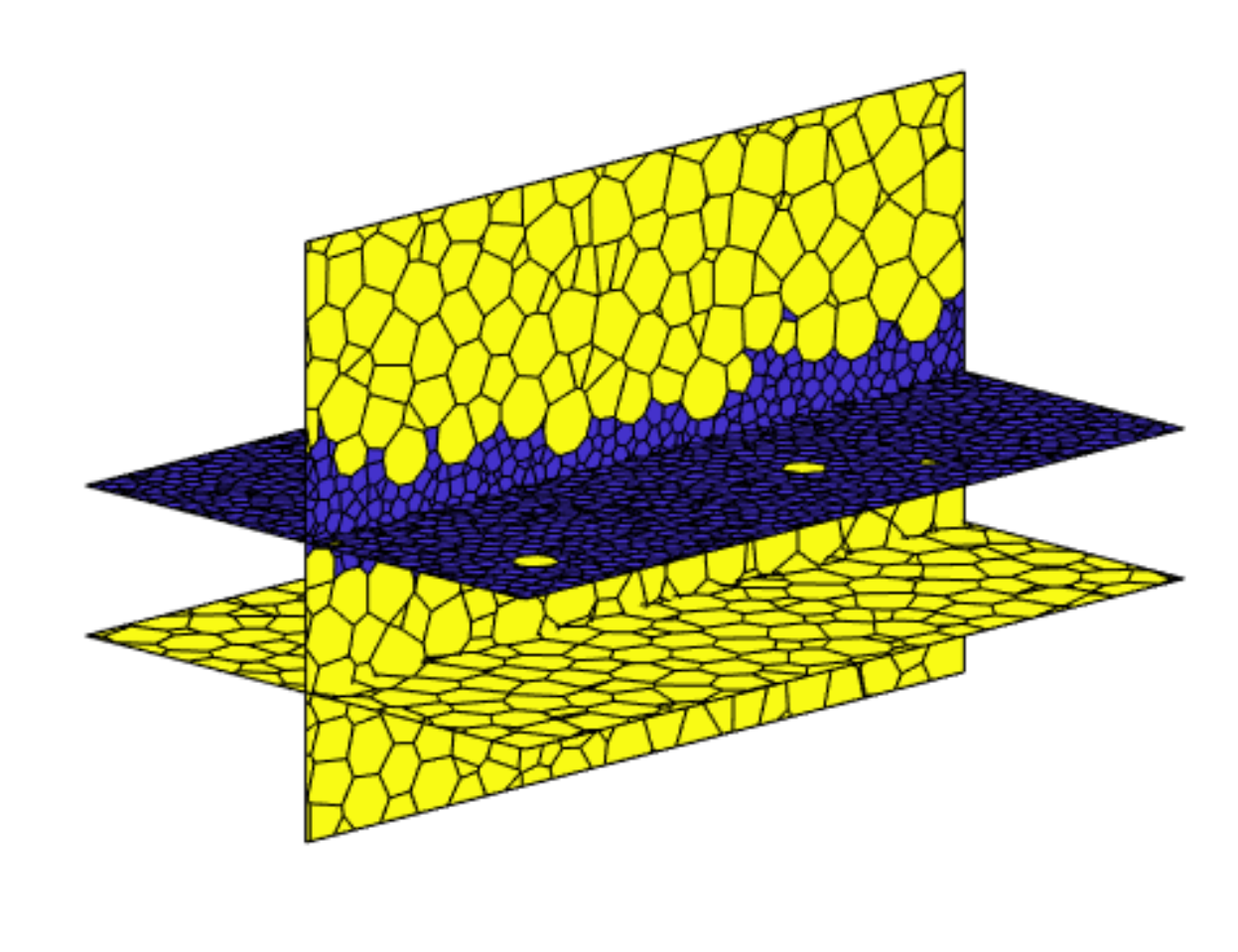}
\caption{\label{Fig:Bands} A periodic RVE of a dual phase material with a banded microstructure (see Example \ref{Example:Bands}).}
\end{figure}

\end{example}

\change{
\begin{example}[Generating an RVE with a log-normal distribution of grain volumes]
	\label{Example:LogNormal}
	Figure \ref{Fig:LogNormal} shows an example of a periodic Laguerre diagram generated by Algorithm 2, in which the grains have volumes that are log-normally distributed. There are $n=10,000$ grains. We used $K=5$ regularisation steps and a volume tolerance of $1\%$. The grains are coloured by their volume, using a log scale. We placed the initial seeds $\bx_1^{(0)},\ldots,\bx_n^{(0)}$ at random. The target volumes were generated by drawing $10,000$ samples $r_i$ from a log-normal distribution with mean $1$ and standard derivation $0.35$ (these correspond to the log-normal parameters $\sigma=(\log(1+0.35^2))^{1/2}=0.34$ and $\mu=-\sigma^2/2$).
	The target volumes were defined by
	\begin{equation}
	m_i=\frac{r_i^3}{L_1L_2L_3\sum_{j=1}^n r_j^3},
\end{equation}
	where the $L_i$ are the side-lengths of the domain $\Omega=[0,L_1]\times[0,L_2]\times[0,L_3]$. For large $n$ these target volumes are approximately log-normally distributed with coefficient of variation (the ratio of standard variation to mean) of
	\begin{equation}
	\frac{1}{3}\sqrt{\exp\left(9\sigma^2\right)-1}.
\end{equation}
	The algorithm took $669$ seconds using the same machine as above. Observe from Example \ref{Example:RunTimeTests} that Algorithm 2 took only 147 seconds to produce a monodisperse RVE and 207 seconds to produce a dual phase RVE with the same number of grains ($n=10,000$) and the same number of regularisation steps ($K=5$). Therefore the run time of Algorithm 2 increases as the RVE becomes more polydisperse.
	\begin{figure}[!htbp]
		\centering
		\includegraphics[trim={2.5cm 1.5cm 2.5cm 0}, clip,width=0.49\textwidth]{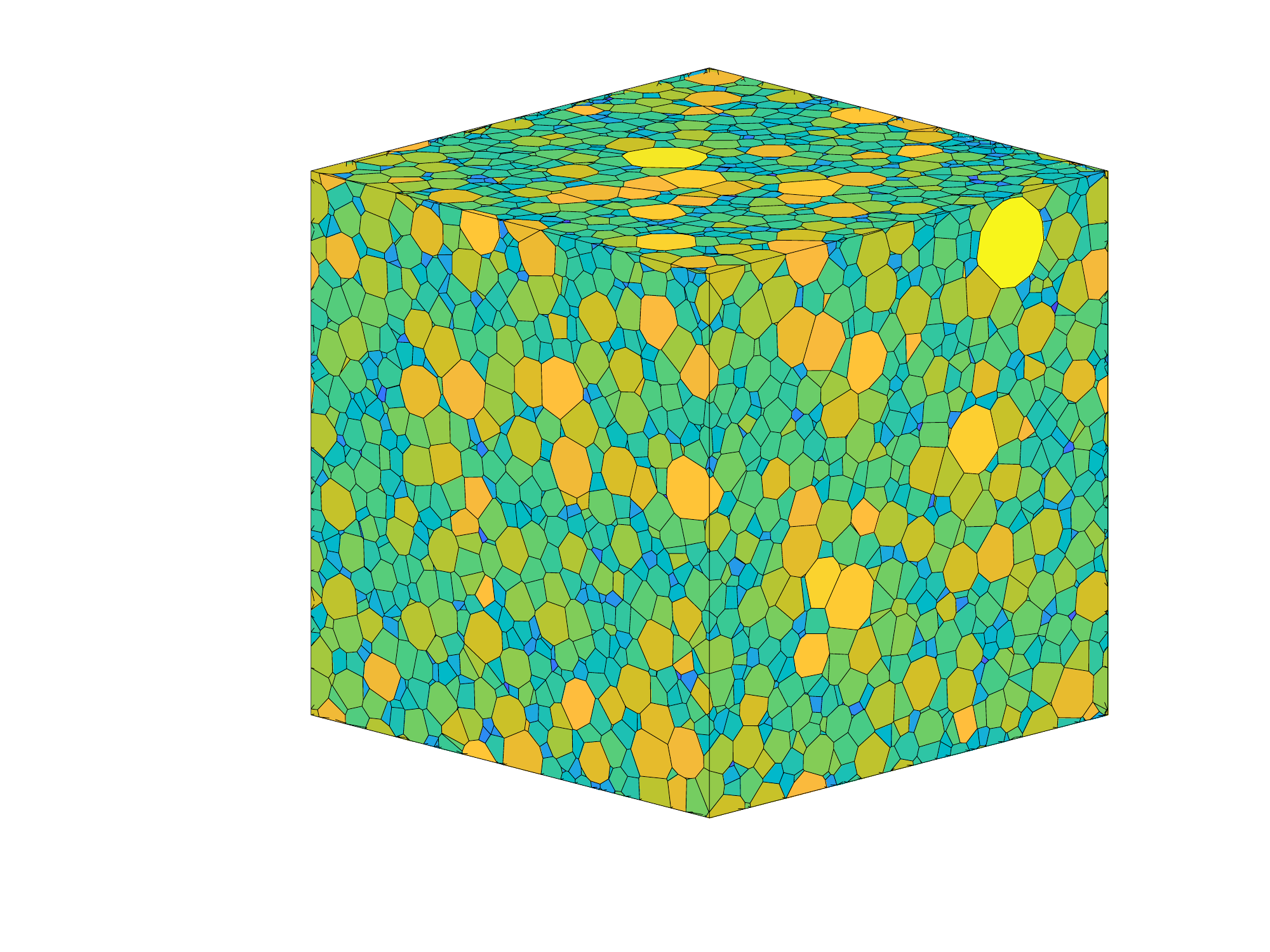}
                \includegraphics[trim={2.5cm 1.5cm 2.5cm 0}, clip,width=0.49\textwidth]{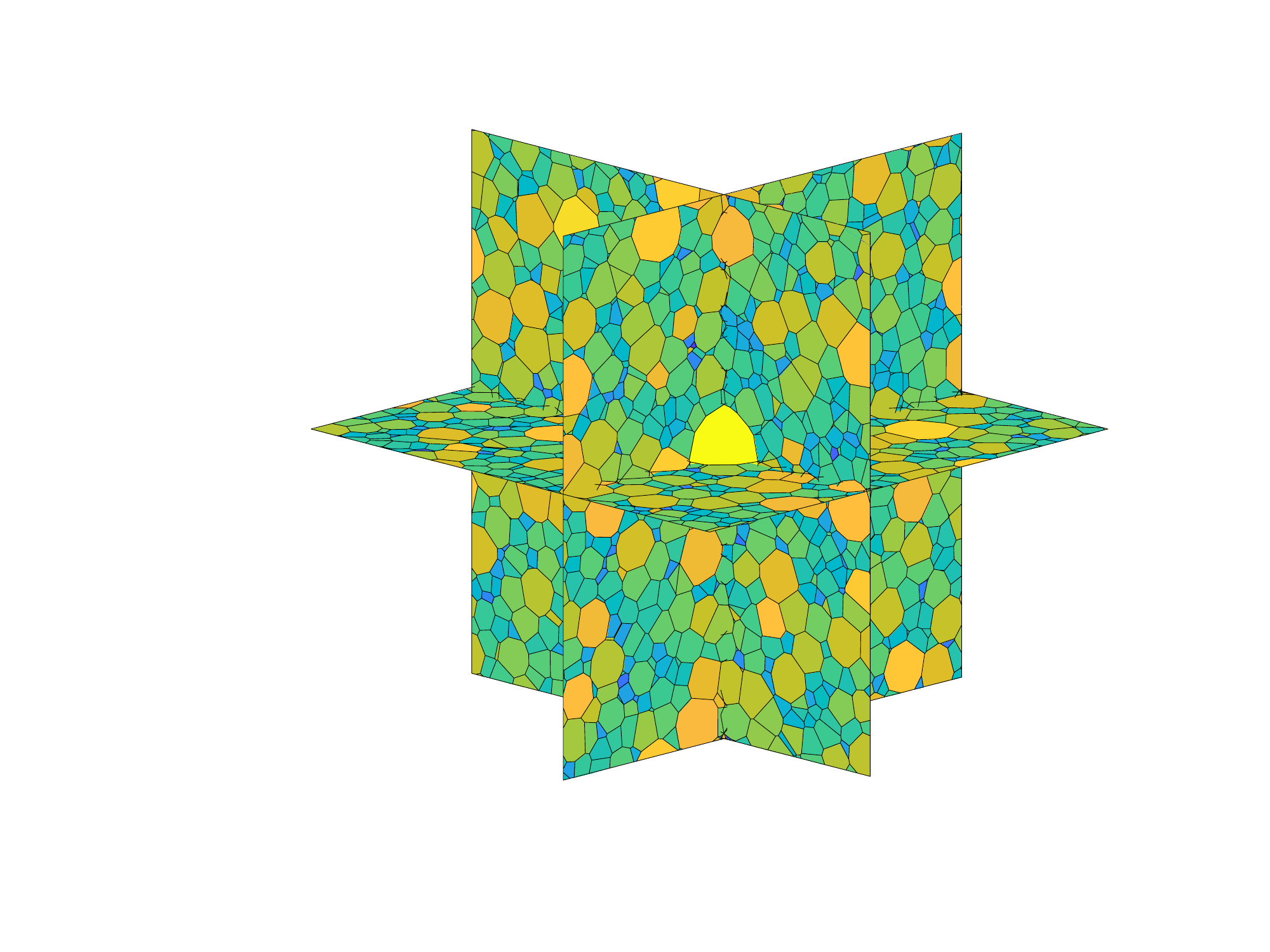}
		\caption{\label{Fig:LogNormal} \change{An RVE in which the grain volumes have approximate\rev{ly} log-normal distribution (see Example \ref{Example:LogNormal}). There are $n=10,000$ grains in the cubic cell. The coefficient of variation of the volumes (the ratio of the standard deviation to the mean) is $1.4$. Cross-sections of the cell are also shown, showing the distribution of sizes throughout the cell.}}
	\end{figure}
\end{example}
}

\section{Discussion}
\change{
	The advantages of our method are
	\begin{itemize}
		\item it is fast;
		\item it can create Laguerre diagrams with grains of \emph{exact} volumes, in principle of any desired tolerance up to machine precision;
		\item it gives some control over the spatial distribution of the grains;
		\item it can create periodic and non-periodic Laguerre diagrams.
	\end{itemize}
	The disadvantages of our method are
	\begin{itemize}
		\item it provides no direct control over the centroids of the grains or their morphology, e.g., their aspect ratio;
		\item it is currently limited to Laguerre diagrams and so the grains cannot have curved boundaries or be non-convex.
	\end{itemize}
	We now discuss these advantages and disadvantages in more detail,
	give evidence in support of our claims, and compare our method with others in the literature.
}

\subsection{Controlling grain volumes: Speed and accuracy}
\change{
	Figure \ref{Fig:RunTimes} \change{shows} that we can create Laguerre diagrams in 3D with up to 20,000 grains in around 10 minutes on a standard desktop PC (\change{without} using parallel computation), where the volumes of the grains are correct to within $1\%$. For 10,000 grains we require as little as 2.5 minutes (see Figure \ref{Fig:RunTimes}), although the time depends on the regularity of the microstructure and whether the material in monodisperse or polydisperse; in Example \ref{Ex:DataFitting} it took 6.25 minutes for 9211 grains and in Example \ref{Example:LogNormal} it took about 11 minutes for 10,000 grains.
}

\change{In our implementation of Algorithm 2 we simply used MATLAB's all-purpose \emph{fminunc} optimisation algorithm. With modern, customised optimal transport optimisation algorithms such as \cite{KitagawaMerigotThibert2019,LevySemiDiscrete2015} it should be possible to use our method to generate Laguerre diagrams with given volume distributions with 100,000 grains in a few minutes \cite[Table 3]{LevySemiDiscrete2015} or even one million grains in less than an hour \cite[Table 4]{LevySemiDiscrete2015}.
}

\change{
	We now compare this with the speed of other algorithms. It is difficult to make a direct comparison in some cases since different methods fit different geometric properties.}

\change{
	In \cite{PetrichEtAl2019} a stochastic optimisation method (the cross-entropy method) is used to solve a non-convex optimisation problem to fit a Laguerre diagram to 3D XRD measurements of grain volumes \emph{and centroids}. \change{The authors} report simulation times (performed using parallel computation) of 19.2 hours for 1439 grains and 122.3 hours for 8063 grains. Note that it is hard to compare their run times with ours since they are also fitting centroids; their method does not try to fit the volumes exactly like we do, but rather obtain a good fit for both volumes and centroids, and
	their method can produce empty Laguerre cells (grains with volume zero).
	While the main focus of our paper is to fit volumes only, we gave an example of fitting volumes and centroids in Example \ref{Ex:DataFitting}, where we fit a Laguerre diagram to 3D EBSD measurements of an IF steel with 9211 grains. The run time is 6.25 minutes and the volumes are correct to within $0.56\%$.
	The centroid errors of most of the grains are comparable to those given in \cite[Fig.~9]{PetrichEtAl2019}, although overall our method does worse than \cite{PetrichEtAl2019} at centroid fitting, as expected.
}

\change{
	Sphere-packing methods are popular for fitting Laguerre diagrams to volume distributions
	\cite{DepriesterKubler2019,LLLFP2011,PFCRMMO2019,WuZhouWangYang2010}. \change{For} $n$ non-overlapping spheres $S_1,\ldots,S_n$ with centres $\bx_1,\ldots,\bx_n$ and radii $r_1,\ldots,r_n$, the Laguerre diagram with seeds $\bx_i$ and weights $w_i=r_i^2$ has the property that cell $L_i$ contains sphere $S_i$. Therefore the volume of $L_i$ is at least the volume of $S_i$. The idea of sphere-packing methods is that if the spheres are close-packed, then the volumes of the Laguerre cells are approximately equal to the volumes of the solid spheres. The disadvantage of this method is that it is inexact and computationally expensive since
	the sphere-packing problem is NP hard \cite{HifiHallah2009}. Nevertheless, this method provided us with inspiration for a good initial guess for the optimisation simulation in Example \ref{Ex:DataFitting} (see also Section \ref{Subsec:InitW}).
}

\change{
	In \cite{AlpersEtAl2015} a method is proposed for fitting grain measurements with \emph{generalised balanced power diagrams} (GBPDs), which are a generalization of Laguerre diagrams. GBPDs are generated by triples $(\bx_i,w_i,A_i)$ of seeds $\bx_i$, weights $w_i$, and positive definite matrices $A_i$; the matrices $A_i$ give some control over the aspect ratio of the generalised Laguerre cells; the case $A_i=I$ for all $i$ corresponds to a standard Laguerre diagram. The advantage of GBPDs is that they give a high degree of control over the morphology of the grains \cite[Figs.~1-6]{AlpersEtAl2015}. The disadvantage is that they are hard to compute. In \cite{AlpersEtAl2015} \change{the authors} approximate GBPDs by voxels, and they fit discretised GBPDs to grain measurements by solving a high-dimensional linear programming problem, where the number of unknowns equals the number of grains multiplied by the number of voxels. \change{It is reported} that to fit a discretised GBPD to 109 grains in 3D took around 6 hours on a standard laptop (this involved solving a linear programming problem with over 77 million variables and 78 million constraints) \cite[Sec.~5.3]{AlpersEtAl2015}.
	\change{Again, it is not possible to make a direct comparison of these run times with the ones presented here since grain volumes and morphology are fitted in \cite{AlpersEtAl2015}, not only grain volumes like here.}
}

\change{
	A heuristic method for approximately fitting GBPDs to measurements of grain volumes, centroids and aspect ratios was proposed \change{in} \cite{TeferraRowenhorst}.
	Their method entirely avoids solving an optimisation problem; it includes explicit formulas for the generators $(\bx_i,w_i,A_i)$ in terms of the data. It is reported in \cite{TeferraRowenhorst} that the method is comparable in accuracy with the optimisation methods of \cite{AlpersEtAl2015,SBDWKKS2016,SBWPBSJ2016} but takes a small fraction of the computation time. No run times or volume errors are reported in \cite{TeferraRowenhorst} \change{precluding} a more precise comparison with our method. \change{Like the sphere-packing method, this heuristic method could be used for generating good initial guesses for optimisation methods.}
}

\subsection{Controlling grain geometry}
\change{
	The main goal of our method is to fit grain volumes \change{quickly} and accurately. Unlike other methods \cite{AlpersEtAl2015,PetrichEtAl2019,SBWPBSJ2016,SDKSJ2017,SWKKS2018,SBDWKKS2016,TeferraRowenhorst} it is not specifically designed for fitting grain morphology. We now discuss to what extent we can control the geometry of Laguerre diagrams.
}

\change{
	Our method gives some control over the spatial distribution of the grains, as shown in Figures \ref{fig:init_seeds_bimodal}, \ref{fig:init_seeds_gradient} \change{and} \ref{Fig:Bands}, where we create microstructures with bands, clusters, and size gradients.
}

\change{
	Several authors use grain centroids as a measure of fitting-error when fitting Laguerre diagrams to data measurements, e.g., \cite{PetrichEtAl2019,TeferraRowenhorst}. We show how we can approximately fit grain centroids to 3D EBSD data in Example \ref{Ex:DataFitting}, although the accuracy is much lower than the volume accuracy.
}

\change{
	In its current form\rev{,} our method gives no direct control over the aspect ratio of the grains. Like the sphere-packing method, Algorithm 2 tends to produce grains that are too round compared to grains typically seen in metals; see Section \ref{Subsec:StopCrit}.
}

\change{
	\rev{Nevertheless} there are several ways \rev{how} our method could be generalised to give more control over the morphology of the grains. For example, by combining our method with
	\emph{multilevel Voronoi diagrams} \cite{KokKorver99,YadegariTurteltaubSuikerKok14} we could maintain control over the volume of the grains while producing more realistic RVEs with non-convex and elongated grains. The idea would be to first use Algorithm 2 to create a Laguerre diagram with $N$ `micro-grains' of equal volume for large $N$.
	Then we would glue together the micro-grains into non-convex `macro-grains'. By choosing which \rev{micro-}grains to glue, we would control the volume and the morphology of the macro-grains. (The multilevel Voronoi method glues together two micro-grains if their generators lie in the same Voronoi cell of a `coarser' Voronoi diagram with fewer generators.)}

\change{
	In principle our algorithms can also be generalised very easily to produce GBPDs with grains of given volumes (up to any desired tolerance) by modifying the objective functions $g$ and $g_k$ in Algorithms 1 and 2 (simply replace the Laguerre cells $L_i$ with generalised Laguerre cells, and replace the isotropic distances $|\bx-\bx_i|$ with anisotropic distances $|\bx-\bx_i|_{A_i}$). This would again allow us to control both the volumes and the aspect ratio of the grains. In practice, however, it is expensive to compute GBPDs to high accuracy; discretizing them with voxels greatly increases the cost of the algorithm. Without developing new computational geometry algorithms for the efficient computation of GBPDs, this limits the method to a small number of grains or greatly increases the run time (cf.~the run time of 6 hours for 109 grains in 3D in \cite{AlpersEtAl2015}).
}

\change{
	Since our method is currently limited to Laguerre diagrams, the grains cannot have curved boundaries or be non-convex. Curved grain boundaries can be created using additively-weighted Voronoi diagrams (Apollonius diagrams) \cite{AurenhammerKleinLee13}, anisotropic diagrams \cite{AltendorfEtAl2014}, or GBPDs \cite{AlpersEtAl2015,SBWPBSJ2016,SWKKS2018,TeferraRowenhorst}, although these are all more costly to compute than Laguerre diagrams. Algorithms 1 and 2 can also be modified to produce Apollonius diagrams with grains of given volumes (in the definition of the objective functions $g$ and $g_k$ simply replace the Laguerre cells $L_i$ with Apollonius cells, and replace the squared distances $|\bx-\bx_i|^2$ with non-squared distances $|\bx-\bx_i|$) but again the implementation cost is an obstacle at the present time. We plan to explore this and the above generalizations in a future paper.}

\section{Conclusions}
\change{
In this paper we introduced a fast optimisation method for generating Laguerre diagrams with user-defined grain size distributions. The volumes of the grains can be controlled exactly (to within any given tolerance). We produced industrially-relevant examples of RVEs with up to 20,000 grains with only $1\%$ volume error in the order of minutes on a standard desktop PC.
We also demonstrated how the spatial and texture distribution of the grains can be partially controlled.
Our algorithms can create both non-periodic Laguerre diagrams (for data fitting) and \rev{periodic} Laguerre diagrams (for generating RVEs of polycrystalline metals or solid foams for computational homogenisation).
}	


\section*{Acknowledgments}
  The authors would like to thank Carola Celada-Casero for useful \change{discussions}. DPB would like to thank the EPSRC for financial support via the grant EP/R013527/1, EP/R013527/2 Designer Microstructure via Optimal Transport Theory. Some of the work of DPB was carried out at Durham University.
The work on generating 3D EBSD data has received funding from the European Union's Horizon 2020 research and innovation programme \emph{Euratom research and training programme 2014-2018} under grant agreement No 709418 MuSTMeF.

\appendix

	\section{Proof that Algorithm 2 is energy-decreasing and convergent}
	\label{Appendix1}
	\change{
		Throughout this section we assume that $\Omega \subset \mathbb{R}^d$ is compact.}
	
	\change{
	First we prove that Algorithm 2 is energy-decreasing, equation \eqref{eq:EnergyDecreasing}. Recall that if $U$ is a compact subset of $\mathbb{R}^d$ with centroid $\bc(U)$, then
	\begin{equation}
	\label{eq:minf}
	\min_{\bz \in \mathbb{R}^d} \int_U |\bx-\bz|^2 \, d \bx =
	\int_U |\bx-\bc(U)|^2 \, d \bx.
	\end{equation}
	This follows from the fact that the function $\bz \mapsto \int_U |\bx-\bz|^2 \, d \bx$ is strictly convex with unique critical point $\bc(U)$.
	By equation \eqref{eq:E} and by the way
	$\big\{ \big( \bx_i^{(k)},w_i^{(k)} \big) \big\}_{i=1}^n$
	is constructed using Algorithm 2,
	\begin{equation}
	E(\bx_1^{(k)},\ldots,\bx_n^{(k)})=
	\sum_{i=1}^n \int_{L_i^{(k)}} \big| \bx - \bx^{(k)}_i \big|^2 \, \d \bx
\end{equation}
	where $\big\{ L_i^{(k)} \big\}_{i=1}^n$ is the Laguerre diagram with generators $\big\{ \big(\bx^{(k)}_i,w^{(k)}_i \big) \big\}_{i=1}^n$.
	Therefore
	\begin{align}
	\nonumber
	E(\bx_1^{(k)},\ldots,\bx_n^{(k)})
	& \stackrel{\eqref{eq:minf}}{\ge} \sum_{i=1}^n \int_{L_i^{(k)}} \big| \bx - \bc\big(L_i^{(k)}\big) \big|^2 \, \d \bx
	\\
	\nonumber
	& \stackrel{\phantom{(3.3)}}{=} \sum_{i=1}^n \int_{L_i^{(k)}} \big| \bx - \bx^{(k+1)}_i \big|^2 \, \d \bx
	& \quad \textrm{(regularisation step of Alg.~2)}
	\\
	& \stackrel{\eqref{eq:E0}}{\ge} E(\bx_1^{(k+1)},\ldots,\bx_n^{(k+1)}).
	\end{align}
	This proves \eqref{eq:EnergyDecreasing}. The inequalities above are strict unless $\bx_i^{(k)}=\bx_i^{(k+1)}$ for all $i$, which means that $\big(\bx_1^{(k)},\ldots,\bx_n^{(k)}\big)$ is a fixed point of Algorithm 2.
}

\change{
	Next we prove a \emph{weak global convergence} result of the form \cite[Theorem 3.8]{EmelianenkoJuRand2008}, where convergence of the classical Lloyd algorithm was proved.
	Weak global convergence means that $\nabla E\big(\bx_1^{(k)},\ldots,\bx_n^{(k)}\big) \to \mathbf{0}$  as $k \to \infty$ and that any convergent subsequence of $\big(\bx_1^{(k)},\ldots,\bx_n^{(k)}\big)$ converges to a critical point of $E$, namely to a centroidal Laguerre diagram. This convergence is called \emph{weak} because it does not guarantee convergence of the whole sequence $\big(\bx_1^{(k)},\ldots,\bx_n^{(k)}\big)$ (different subsequences could converge to different critical points).
}

\change{
	By construction $\bx_i^{(k)}$ is the centroid of the convex set $L_i^{(k-1)}$, which has volume $m_i$. Therefore by
	\cite[Lemma 3.2]{EmelianenkoJuRand2008} the distance between
	$\bx_i^{(k)}$ and $\partial L_i^{(k-1)}$ \rev{has a lower bound of} 
	$C m_i^2/\mathrm{diam}(\Omega)^{2d-1}$ where $C=1/1024$.
	Therefore the closest that two generators $x_i^{(k)}$ and $x_j^{(k)}$ can be is  $2 C m^2/\mathrm{diam}(\Omega)^{2d-1}$ where $m=\min_i m_i$. Note that this bound is independent of the iteration number $k$. Therefore the iterates $(\bx_1^{(k)},\ldots,\bx_n^{(k)})$ lie in the compact set
	\begin{equation}
	\label{eq:compactset}
	\{ (\bx_1,\ldots,\bx_n) \in \Omega^n : |\bx_i - \bx_j| \ge 2 C m^2/\mathrm{diam}(\Omega)^{2d-1} \; \forall \, i,j\}.
	\end{equation}
	\rev{Owing} to this compactness and the energy-decreasing property of Algorithm 2, we have weak global convergence of  $\big(\bx_1^{(k)},\ldots,\bx_n^{(k)}\big)$; see the Global Convergence Theorem in \cite[p.~206]{LuenbergerYe} or \cite[proof of Theorem 3.3]{BourneRoper15}.
}

\change{
	Finally we prove a \emph{strong convergence} result, namely that the whole sequence
	$\big(\bx_1^{(k)},\ldots,\bx_n^{(k)}\big)$ produced by Algorithm 2 converges to a critical point of $E$. We are only able to prove this\rev{,} however\rev{,} under the following generic assumption:
	There are only finitely many centroidal Laguerre diagrams with the same energy $E$. More precisely we assume that, for all $M>0$,
	\begin{equation}
	\label{eq:ass}
	\# \left\{ (\bx_1,\ldots,\bx_n) \in \Omega^n : \nabla E(\bx_1,\ldots,\bx_n) = \mathbf{0}, \; E(\bx_1,\ldots,\bx_n) = M \right\} < \infty.
	\end{equation}
	This assumption is expected to hold for `generic' domains $\Omega$ and masses $m_1,\ldots,m_n$ \cite[p.~107]{DuEmelianenkoJu},
	\cite[Remark 3.4]{BourneRoper15}, however there are examples where it fails. For example, if $\Omega$ is a disc, $n=2$ and $m_1=m_2$, then there are infinitely many critical points of $E$ with the same energy by the rotational symmetry of $\Omega$. Note that if $E$ satisfies the generic condition of being a Morse function (having no degenerate critical points), then its critical points are isolated. Since they lie in the compact set \eqref{eq:compactset} there are only finitely many of them, and so assumption \eqref{eq:ass} is satisfied.
}

\change{
	The Monotone Convergence Theorem implies that the whole sequence
	$E\big(\bx_1^{(k)},\ldots,\bx_n^{(k)}\big)$ converges (because Algorithm 2 is energy-decreasing and $E$ is bounded below by zero). In addition $E$ is continuous, and so
	every accumulation point of
	$\big(\bx_1^{(k)},\ldots,\bx_n^{(k)}\big)$ has the same energy $E$. Moreover, by the global weak convergence result above, every  accumulation point is a critical point of $E$. Therefore  assumption \eqref{eq:ass} ensures there are only finitely many accumulation points.}

\change{	
	We now complete the proof following the idea from \cite[proof of Theorem 2.5]{DuEmelianenkoJu}.
	Assume for contradiction that the sequence  $\bX^{(k)}:=\big(\bx_1^{(k)},\ldots,\bx_n^{(k)}\big)$ does not converge. Since it only has finitely many accumulation points, there exist
	distinct accumulation points  $\bY:=(\by_1,\ldots,\by_n)$, $\bZ:=(\bz_1,\ldots,\bz_n)$ and distinct subsequences $\bX^{(k_j)}:=\big(\bx_1^{(k_j)},\ldots,\bx_n^{(k_j)}\big)$, $\bX^{(k_l)}:=\big(\bx_1^{(k_l)},\ldots,\bx_n^{(k_l)}\big)$ such that $\bX^{(k_j)} \to \bY$, $\bX^{(k_l)}\to \bZ$, and $k_l=k_j+1$, i.e., the sequence $\bX^{(k)}$ `jumps' between the two subsequences infinitely many times. (Note that such subsequences may not exist if there are infinitely many accumulation points.) Let $\delta = |\bY-\bZ|>0$ and let $T:\Omega^n \to \Omega^n$ denote the continuous map that sends
	$(\bx_1,\ldots,\bx_n)$ to $(\bc(L_1),\ldots,\bc(L_n))$, where $\{ L_i\}$ is the Laguerre diagram with seeds $(\bx_1,\ldots,\bx_n)$ and cells of volume $m_1,\ldots,m_n$. In other words, $\bX^{(k+1)}=T(\bX^{(k)})$. Moreover, $T(\bY)=\bY$ and $T(\bZ)=\bZ$. Then
	\begin{align}
	\nonumber
	\delta & = |\bY - \bZ|
	\\
	\nonumber
	& \le
	| \bY -\bX^{(k_j)}| + |\bX^{(k_j)}-\bX^{(k_l)}| + | \bX^{(k_l)}-\bZ|
	\\
	& = | \bY -\bX^{(k_j)}| + |\bX^{(k_j)}-T(\bX^{(k_j)})| + | \bX^{(k_l)}-\bZ|.
	\end{align}
	By taking $j,l \to \infty$ in the right-hand side, and using the continuity of $T$, we find that $\delta=0$, which is a contradiction.
	Therefore the whole sequence  $\big(\bx_1^{(k)},\ldots,\bx_n^{(k)}\big)$ converges to a critical point of $E$. Moreover, this critical point must a local minimum point or a saddle point of $E$ by the energy-decreasing property. This completes the proof.
}

\bibliographystyle{tfq}
\bibliography{RVEpaper_bib}

\begin{thebibliography}{10}
\newcommand{\printfirst}[2]{#1}
\newcommand{\switchargs}[2]{#2#1}
\providecommand{\url}[1]{\normalfont{#1}}
\providecommand{\urlprefix}{Available at }

\bibitem{AlpersEtAl2015}
A. Alpers, A. Brieden, P. Gritzmann, A. Lyckegaard, and H.F. Poulsen,
  \emph{Generalized balanced power diagrams for 3{D} representations of
  polycrystals}, Philosophical Magazine 95 (2015), pp. 1016--1028,
  \urlprefix\url{https://doi.org/10.1080/14786435.2015.1015469}.

\bibitem{AltendorfEtAl2014}
H. Altendorf, F. Latourte, D. Jeulin, M. Faessel, and L. Saintoyant, \emph{3{D}
  reconstruction of a multiscale microstructure by anisotropic tessellation
  models}, Image Analysis \& Stereology 33 (2014), pp. 121--130,
  \urlprefix\url{https://www.ias-iss.org/ojs/IAS/article/view/1090}.

\bibitem{BarkerBollerheyHamaekers}
J. Barker, G. Bollerhey, and J. Hamaekers, \emph{A multilevel approach to the
  evolutionary generation of polycrystalline structures}, Computational
  Materials Science 114 (2016), pp. 54--63,
  \urlprefix\url{http://www.sciencedirect.com/science/article/pii/S0927025615007259}.

\bibitem{DepriesterKubler2019}
D. Depriester and R. Kubler, \emph{Radical {V}oronoi tessellation from random
  pack of polydisperse spheres: {P}rediction of the cells' size distribution},
  Computer-Aided Design 107 (2019), pp. 37--49,
  \urlprefix\url{https://doi.org/10.1016/j.cad.2018.09.001}.

\bibitem{KokKorver99}
P.J.J. Kok and F.N.M. Korver, \emph{Modelling of complex microstructures in
  multi phase steels: geometric considerations for building an RVE}, in
  \emph{Proceedings of the {X} International Conference on Computational
  Plasticity}. 2009.

\bibitem{Liebscher2015}
A. Liebscher, \emph{Laguerre approximation of random foams}, Philosophical
  Magazine 95 (2015), pp. 2777--2792,
  \urlprefix\url{https://doi.org/10.1080/14786435.2015.1078511}.

\bibitem{LeonardiScardiLeoni}
A. Leonardi, P. Scardi, and M. Leoni, \emph{Realistic nano-polycrystalline
  microstructures: beyond the classical {V}oronoi tessellation}, Philosophical
  Magazine 92 (2012), pp. 986--1005,
  \urlprefix\url{https://doi.org/10.1080/14786435.2011.637984}.

\bibitem{LLLFP2011}
A. Lyckegaard, E.M. Lauridsen, W. Ludwig, R.W. Fonda, and H.F. Poulsen,
  \emph{On the use of {L}aguerre tessellations for representations of 3d grain
  structures}, Advanced Engineering Materials 13 (2011), pp. 165--170,
  \urlprefix\url{https://onlinelibrary.wiley.com/doi/abs/10.1002/adem.201000258}.

\bibitem{PetrichEtAl2019}
L. Petrich, J. Stan\v{e}k, M. Wang, D. Westhoff, L. Heller, P. \v{S}ittner,
  C.E. {Krill III}, V. Bene\v{s}, and V. Schmidt, \emph{Reconstruction of
  grains in polycrystalline materials from incomplete data using {L}aguerre
  tessellations}, Microscopy and Microanalysis 25 (2019), pp. 743--752,
  \urlprefix\url{https://doi.org/10.1017/S1431927619000485}.

\bibitem{PFCRMMO2019}
I. P\'erez, M. Muniz~de  Farias, M. Castro, R. Rosell\'o, C. Recarey~Morfa, L.
  Medina, and E. O\~nate, \emph{Modeling polycrystalline materials with
  elongated grains}, International Journal for Numerical Methods in Engineering
  118 (2019), pp. 121--131,
  \urlprefix\url{https://onlinelibrary.wiley.com/doi/abs/10.1002/nme.6004}.

\bibitem{SBDWKKS2016}
A. Spettl, T. Brereton, Q. Duan, T. Werz, C.E. {Krill III}, D.P. Kroese, and V.
  Schmidt, \emph{Fitting {L}aguerre tessellation approximations to tomographic
  image data}, Philosophical Magazine 96 (2016), pp. 166--189,
  \urlprefix\url{https://doi.org/10.1080/14786435.2015.1125540}.

\bibitem{SWKKS2018}
O. \v{S}ediv\'{y}, D. Westhoff, J. Kope\v{c}ek, C.E. Krill III, and V. Schmidt,
  \emph{Data-driven selection of tessellation models describing polycrystalline
  microstructures}, Journal of Statistical Physics 172 (2018), pp. 1223--1246,
  \urlprefix\url{https://doi.org/10.1007/s10955-018-2096-8}.

\bibitem{TeferraRowenhorst}
K. Teferra and D.J. Rowenhorst, \emph{Direct parameter estimation for
  generalised balanced power diagrams}, Philosophical Magazine Letters 98
  (2018), pp. 79--87,
  \urlprefix\url{https://doi.org/10.1080/09500839.2018.1472399}.

\bibitem{WuCaoFan2005}
Y. Wu, J. Cao, and Z. Fan, \emph{Chord length distribution of {V}oronoi diagram
  in {L}aguerre geometry with lognormal-like volume distribution}, Materials
  Characterization 55 (2005), pp. 332--339,
  \urlprefix\url{http://www.sciencedirect.com/science/article/pii/S1044580305001695}.

\bibitem{AlsayednoorHarrisonGuo13}
J. Alsayednoor, P. Harrison, and Z. Guo, \emph{Large strain compressive
  response of 2-{D} periodic representative volume element for random foam
  microstructures}, Mechanics of Materials 66 (2013), pp. 7--20,
  \urlprefix\url{http://www.sciencedirect.com/science/article/pii/S016766361300118X}.

\bibitem{GhoshDimiduk}
S. Ghosh and D. Dimiduk, \emph{Computational Methods for
  Microstructure-Property Relationships}, Springer, 2011.

\bibitem{YadegariTurteltaubSuikerKok14}
S. Yadegari, S. Turteltaub, A.S.J. Suiker, and P.J.J. Kok, \emph{Analysis of
  banded microstructures in multiphase steels assisted by
  transformation-induced plasticity}, Computational Materials Science 84
  (2014), pp. 339--349,
  \urlprefix\url{http://www.sciencedirect.com/science/article/pii/S0927025613007519}.

\bibitem{WuZhouWangYang2010}
Y. Wu, W. Zhou, B. Wang, and F. Yang, \emph{Modeling and characterization of
  two-phase composites by {V}oronoi diagram in the {L}aguerre geometry based on
  random close packing of spheres}, Computational Materials Science 47 (2010),
  pp. 951--961,
  \urlprefix\url{http://www.sciencedirect.com/science/article/pii/S092702560900439X}.

\bibitem{DuanEtAl}
Q. Duan, D. Kroese, T. Brereton, A. Spettl, and V. Schmidt, \emph{Inverting
  {L}aguerre tessellations}, The Computer Journal 57 (2014), pp. 1431--1440.

\bibitem{Teferra2015}
K. Teferra and L. Graham-Brady, \emph{Tessellation growth models for
  polycrystalline microstructures}, Computational Materials Science 102 (2015),
  pp. 57--67, \urlprefix\url{https://doi.org/10.1016/j.commatsci.2015.02.006}.

\bibitem{SBWPBSJ2016}
O. \v{S}ediv\'{y}, T. Brereton, D. Westhoff, L. Pol\'{i}vka, V. Bene\v{s}, V.
  Schmidt, and A. J\"{a}ger, \emph{3{D} reconstruction of grains in
  polycrystalline materials using a tessellation model with curved grain
  boundaries}, Philosophical Magazine 96 (2016), pp. 1926--1949,
  \urlprefix\url{https://doi.org/10.1080/14786435.2016.1183829}.

\bibitem{Dream3D}
M. Groeber and M. Jackson, \emph{{DREAM.3D}: {A} digital representation
  environment for the analysis of microstructure in 3{D}}, Integrating
  Materials 3 (2014), pp. 56--72,
  \urlprefix\url{https://doi.org/10.1186/2193-9772-3-5}.

\bibitem{Neper}
R. Quey and L. Renversade, \emph{Optimal polyhedral description of 3d
  polycrystals: Method and application to statistical and synchrotron x-ray
  diffraction data}, Computer Methods in Applied Mechanics and Engineering 330
  (2018), pp. 308--333,
  \urlprefix\url{http://www.sciencedirect.com/science/article/pii/S0045782517307028}.

\bibitem{SDKSJ2017}
O. \v  Sediv\'{y}, J. Dake, C.E. {Krill III}, V. Schmidt, and A. J\"{a}ger,
  \emph{Description of the 3{D} morphology of grain boundaries in aluminum
  alloys using tessellation models generated by ellipsoids}, Image Analysis \&
  Stereology 36 (2017), pp. 5--13,
  \urlprefix\url{https://www.ias-iss.org/ojs/IAS/article/view/1656}.

\bibitem{LevySchwindt2018}
B. L\'evy and E.L. Schwindt, \emph{Notions of optimal transport theory and how
  to implement them on a computer}, Computers \& Graphics 72 (2018), pp.
  135--148,
  \urlprefix\url{http://www.sciencedirect.com/science/article/pii/S0097849318300098}.

\bibitem{Santambrogio}
F. Santambrogio, \emph{Optimal transport for applied mathematicians},
  Birkh\"{a}user/Springer, 2015,
  \urlprefix\url{https://doi.org/10.1007/978-3-319-20828-2}.

\bibitem{AurenhammerKleinLee13}
 {F. Aurenhammer, R. Klein, and {D.-T.} Lee}, \emph{Voronoi diagrams and
  {D}elaunay triangulations}, World Scientific, 2013,
  \urlprefix\url{https://doi.org/10.1142/8685}.

\bibitem{OkabeBootsSugiharaChiu}
A. Okabe, B. Boots, K. Sugihara, and S.N. Chiu, \emph{Spatial tessellations:
  {C}oncepts and applications of {V}oronoi diagrams}, 2nd ed., Wiley, 2000,
  \urlprefix\url{https://doi.org/10.1002/9780470317013}.

\bibitem{Aurenhammer98}
F. Aurenhammer, F. Hoffmann, and B. Aronov, \emph{Minkowski-type theorems and
  least-squares clustering}, Algorithmica 20 (1998), pp. 61--76,
  \urlprefix\url{https://doi.org/10.1007/PL00009187}.

\bibitem{Merigot2011}
Q. M\'erigot, \emph{A multiscale approach to optimal transport}, Computer
  Graphics Forum 30 (2011), pp. 1583--1592,
  \urlprefix\url{https://onlinelibrary.wiley.com/doi/abs/10.1111/j.1467-8659.2011.02032.x}.

\bibitem{BoydVandenberghe}
S. Boyd and L. Vandenberghe, \emph{Convex optimization}, Cambridge University
  Press, 2004, \urlprefix\url{https://doi.org/10.1017/CBO9780511804441}.

\bibitem{LevySemiDiscrete2015}
B. L\'{e}vy, \emph{A numerical algorithm for {$L_2$} semi-discrete optimal
  transport in 3{D}}, ESAIM. Mathematical Modelling and Numerical Analysis 49
  (2015), pp. 1693--1715, \urlprefix\url{https://doi.org/10.1051/m2an/2015055}.

\bibitem{KitagawaMerigotThibert2019}
J. Kitagawa, Q. M\'{e}rigot, and B. Thibert, \emph{Convergence of a {N}ewton
  algorithm for semi-discrete optimal transport}, Journal of the European
  Mathematical Society (JEMS) 21 (2019), pp. 2603--2651,
  \urlprefix\url{https://doi.org/10.4171/JEMS/889}.

\bibitem{BriedenGritzmann}
A. Brieden and P. Gritzmann, \emph{On optimal weighted balanced clusterings:
  {G}ravity bodies and power diagrams}, SIAM Journal on Discrete Mathematics 26
  (2012), pp. 415--434, \urlprefix\url{https://doi.org/10.1137/110832707}.

\bibitem{BourneRoper15}
D.P. Bourne and S.M. Roper, \emph{Centroidal power diagrams, {L}loyd's
  algorithm, and applications to optimal location problems}, SIAM Journal on
  Numerical Analysis 53 (2015), pp. 2545--2569,
  \urlprefix\url{https://doi.org/10.1137/141000993}.

\bibitem{LevyCentroidalPower}
 {S.{-}Q. Xin, B. L{\'e}vy, Z. Chen, L. Chu, Y. Yu, C. Tu, and W. Wang},
  \emph{Centroidal power diagrams with capacity constraints: {C}omputation,
  applications, and extension}, ACM Transactions on Graphics 35 (2016), pp.
  244:1--244:12, \urlprefix\url{http://doi.acm.org/10.1145/2980179.2982428}.

\bibitem{DuFaberGunzburger1999}
Q. Du, V. Faber, and M. Gunzburger, \emph{Centroidal {V}oronoi tessellations:
  applications and algorithms}, SIAM Review 41 (1999), pp. 637--676,
  \urlprefix\url{https://doi.org/10.1137/S0036144599352836}.

\bibitem{BollobasStern}
B. Bollob\'{a}s and N. Stern, \emph{The optimal structure of market areas},
  Journal of Economic Theory 4 (1972), pp. 174--179,
  \urlprefix\url{https://doi.org/10.1016/0022-0531(72)90147-0}.

\bibitem{Gruber04}
P.M. Gruber, \emph{Optimum quantization and its applications}, Advances in
  Mathematics 186 (2004), pp. 456--497,
  \urlprefix\url{https://doi.org/10.1016/j.aim.2003.07.017}.

\bibitem{GrayNeuhoff}
R.M. Gray and D.L. Neuhoff, \emph{Quantization}, Institute of Electrical and
  Electronics Engineers. Transactions on Information Theory 44 (1998), pp.
  2325--2383, \urlprefix\url{https://doi.org/10.1109/18.720541}, information
  theory: 1948--1998.

\bibitem{GrafLuschgy2000}
S. Graf and H. Luschgy, \emph{Foundations of quantization for probability
  distributions}, Vol. 1730, Springer, 2000,
  \urlprefix\url{https://doi.org/10.1007/BFb0103945}.

\bibitem{Voro++}
C.H. Rycroft, \emph{{V}{O}{R}{O}++: A three-dimensional {V}oronoi cell library
  in {C}++}, Chaos: An Interdisciplinary Journal of Nonlinear Science 19
  (2009), p. 041111, \urlprefix\url{https://doi.org/10.1063/1.3215722}.

\bibitem{PowerBounded}
 Firman, \emph{Fast bounded power diagram},
  \url{https://uk.mathworks.com/matlabcentral/fileexchange/56633-fast-bounded-power-diagram},
  MATLAB Central File Exchange.

\bibitem{Aurenhammer87}
F. Aurenhammer, \emph{Power diagrams: properties, algorithms and applications},
  SIAM Journal on Computing 16 (1987), pp. 78--96,
  \urlprefix\url{https://doi.org/10.1137/0216006}.

\bibitem{PowerDiagramWrapper}
F. McCollum, \emph{Power diagrams},
  \url{https://uk.mathworks.com/matlabcentral/fileexchange/44385-power-diagrams},
  MATLAB Central File Exchange.

\bibitem{Spettl2014}
A. Spettl, T. Werz, C. {Krill III}, and V. Schmidt, \emph{Parametric
  representation of 3d grain ensembles in polycrystalline microstructures},
  Journal of Statistical Physics 154 (2014), pp. 913--928,
  \urlprefix\url{https://doi.org/10.1007/s10955-013-0893-7}.

\bibitem{HifiHallah2009}
M. Hifi and R. M'Hallah, \emph{A literature review on circle and sphere packing
  problems: Models and methodologies}, Advances in Operations Research 2009
  (2009), pp. 150624:1--150624:22.

\bibitem{EmelianenkoJuRand2008}
M. Emelianenko, L. Ju, and A. Rand, \emph{Nondegeneracy and weak global
  convergence of the {L}loyd algorithm in {$\mathbb{R}^d$}}, SIAM Journal on
  Numerical Analysis 46 (2008), pp. 1423--1441,
  \urlprefix\url{https://doi.org/10.1137/070691334}.

\bibitem{LuenbergerYe}
D. Luenberger and Y. Ye, \emph{Linear and Nonlinear Programming}, 3rd ed.,
  Springer, 2008.

\bibitem{DuEmelianenkoJu}
Q. Du, M. Emelianenko, and L. Ju, \emph{Convergence of the {L}loyd algorithm
  for computing centroidal {V}oronoi tessellations}, SIAM Journal on Numerical
  Analysis 44 (2006), pp. 102--119,
  \urlprefix\url{https://doi.org/10.1137/040617364617364}.

\end{thebibliography}

\end{document}